\documentclass[11pt]{article}

 \usepackage{setspace}
\usepackage{natbib}

\usepackage{latexsym}
\usepackage{calc}

\usepackage{amssymb}
\usepackage{amsmath}
\usepackage{graphicx}

\usepackage{psfrag}

\usepackage{float}

\usepackage{color}

\usepackage{nicefrac}

\usepackage{ntheorem}

\usepackage{hyperref}

\usepackage{tikz}
\usetikzlibrary{backgrounds}

\newcounter{hours}\newcounter{minutes}

\def\nr{\par }

\def\beq{\begin{equation}}
\def\eeq{\end{equation}}

\newtheorem{theorem}{Theorem}

\newtheorem{lemma}{Lemma}
\newtheorem{corollary}{Corollary}

\newtheorem{proposition}{Proposition}
\newtheorem{assumption}{Assumption}
\newtheorem{definition}{Definition}

\newtheorem{example}{Example}
\newtheorem{remark}{Remark}
\newcommand{\proof}{\bf Proof: \rm \nr}

\newcommand{\qed}{\hfill $\Box$ \nr \medskip}

\def\ba{\begin{array}}
\def\ea{\end{array}}
\def\beann{\begin{eqnarray*}}
\def\eeann{\end{eqnarray*}}
\def\bea{\begin{eqnarray}}
\def\eea{\end{eqnarray}}

\def\BT{\begin{theorem}}
\def\ET{\end{theorem}}
\def\BL{\begin{lemma}}
\def\EL{\end{lemma}}
\def\BC{\begin{corollary}}
\def\EC{\end{corollary}}
\def\BE{\begin{example}}
\def\EE{\end{example}}
\def\BD{\begin{definition}}
\def\ED{\end{definition}}
\def\BR{\begin{remark}}
\def\ER{\end{remark}}
\def\BAS{\begin{assumption}}
\def\EAS{\end{assumption}}
\def\BI{\begin{itemize}}
\def\EI{\end{itemize}}

\def\BMP{\begin{minipage}{9.5cm}}
\def\EMP{\end{minipage}}
\def\MPT{\begin{minipage}{11.5cm}}
\def\EPT{\end{minipage}}

\def\H{\mathbb{H}}
\def\R{\mathbb{R}}

\newcommand*\diff{\mathop{}\!\mathrm{d}}

\newcommand*\samethanks[1][\value{footnote}]{\footnotemark[#1]}
\newcommand{\NNorm}[2]{\left\Vert {#1} \right\Vert_{#2}}

\title{Cardinality-constrained optimization problems \\ in general position and beyond}
\author{
S. L\"ammel
\thanks{
Department of Mathematics, Chemnitz University of Technology,
Reichenhainer Str. 41, 09126
Chemnitz, Germany; e-mail: sebastian.laemmel@mathematik.tu-chemnitz.de, vladimir.shikhman@mathematik.tu-chemnitz.de (corresponding author).
 } \and V. Shikhman\samethanks[1]
}

\begin{document}
\thispagestyle{empty}
\maketitle
\vspace{-5ex}
\abstract{
We study cardinality-constrained optimization problems (CCOP) in general position, i.\,e. those optimization-related properties that are fulfilled for a dense and open subset of their defining functions. We show that the well-known cardinality-constrained linear independence constraint qualification (CC-LICQ) is generic in this sense. For M-stationary points we define nondegeneracy and show that it is a generic property too. In particular, the sparsity constraint turns out to be active at all minimizers of a generic CCOP. Moreover, we describe the global structure of CCOP in the sense of Morse theory, emphasizing the strength of the generic approach. Here, we prove that multiple cells need to be attached, each of dimension coinciding with the proposed M-index of nondegenerate M-stationary points.
Beyond this generic viewpoint, we study singularities of CCOP. For that, the relation between nondegeneracy and strong stability in the sense of \cite{kojima:1980} is examined. We show that nondegeneracy implies the latter, while the reverse implication is in general not true. To fill the gap, we fully characterize the strong stability of M-stationary points under CC-LICQ by first- and second-order information of CCOP defining functions. Finally, we compare nondegeneracy and strong stability of M-stationary points with second-order sufficient conditions recently introduced in the literature.
%examine known stationarity concepts and their behaviour in terms of strong stability. 

%, which we define for nondegenerate M-stationary points. The attachment of multiple cells is a special characteristic of CCOP and reflects the difficulty of the solving of CCOP.
}

\vspace{2ex}
{\bf Keywords: general position, nondegenerate M-stationarity, strong stability, second-order sufficient condition, Morse theory}

\vspace{2ex}
{\bf MSC-classification: 90C26, 90C46}

\section{Introduction}
\label{sec:intro}
We consider the class of cardinality-constrained optimization problems:
\[
\mbox{CCOP}: \quad
\min_{x} \,\, f(x)\quad \mbox{s.\,t.} \quad x \in M
\]
with the feasible set given by equality, inequality, and cardinality constraints:
\begin{equation}
    \label{eq:cc2}
    M=\left\{x \in\R^n \left\vert\; h(x)=0, g(x)\ge 0, \left\|x\right\|_0\le s \right.\right\},
\end{equation}
where the so-called $\ell_0$ ”norm” is counting non-zero entries of $x$:
\[
\left\|x\right\|_0 = \left|\left\{i \in \{1,\ldots, n\}\; \vert\; x_i\ne 0\right\}\right|.
\]
Here, we assume that the objective function $f \in C^2(\R^n,\R)$, as well as the equality and inequality constraints 
$h=\left(h_p, p \in P\right) \in C^2(\R^n,\R^{|P|})$,
$g=\left(g_q, q \in Q\right) \in C^2(\R^n,\R^{|Q|})$ are twice continuously differentiable, and $s \in \{0,1,\ldots, n-1\}$ is an integer. We call $f,h,g$ the defining functions and highlight the dependence on them by writing $\mbox{CCOP}(f,g,h)$ if needed. 
%Note that the requirement of sparsity in (\ref{eq:cc2}) is due to various applications, such as compressed sensing, model selection, image processing etc. 

% General position
In this paper, we study CCOP in general position. Loosely speaking, the latter refers to an open and dense subset of its defining functions. In order to state genericity in mathematically precise terms, let $C^2\left(\R^n,\R\right)$ be endowed with the strong (or Whitney) $C^2$-topology, denoted by $C^2_s$, see e.\,g. \cite{hirsch:1976}. The $C^2_s$-topology is generated by allowing perturbations of the functions, their gradients and Hessians, which are controlled by means of continuous positive functions. 
%We say that a set is $C^2_s$-generic if it contains a countable intersection of $C^2_s$-open and -dense subsets. Since $C^2\left(\R^n,\R\right)$ endowed with the $C^2_s$-topology is a Baire space, generic sets are in particular dense. 
Let the product space $$C^2(\R^n, \R^{|P|+|Q|+1}) \cong C^2(\R^n,\R) \times
C^2(\R^n,\R^{|P|}) \times C^2(\R^n,\R^{|Q|})$$ of CCOP defining functions be topologized with the corresponding product topology. We say that a property is generic for CCOP if there exists a $C_s^2$-open and $C_s^2$-dense subset $\mathcal{D} \subset C^2(\R^n,\R) \times
C^2(\R^n,\R^{|P|}) \times C^2(\R^n,\R^{|Q|})$, such that the corresponding $\mbox{CCOP}(f,g,h)$ fulfils this property for all $(f,g,h) \in \mathcal{D}$.
We refer to $\mbox{CCOP}(f,g,h)$ being not in general position as describing singularities. 

The idea behind this distinction is that a unified theory for optimization problems in general position can be established. This includes suitable constraint qualifications, first- and second-order optimality conditions, local and global structure of stationary points, convergence of Newton-type methods, stability w.r.t. data perturbations etc. The classification and description of singularities is more involved. An every singularity often requires the development of an independent theory. This is in view of the fact that singularities are  unstable and tend to change their type due to bifurcations. The main challenge is, thus, to identify generic properties of an optimization problem and to subsequently analyze its singularities, one by one. This generic perspective has been pioneered in \cite{jongen:1977} and popularized in \cite{Jongen:2000} in context of nonlinear programming. Since then, the generic paradigm in optimization has been successively applied to disjunctive programming, general semiinfinite programming, mathematical programs with complementarity constraints, bilevel optimization, mathematical programs with vanishing constraints, generalized Nash equilibrium problems, nonlinear semidefinite programming, mathematical programs with switching constraints etc., see e.\,g. \cite{Shikhman:2012} and references therein.

% Results
Our findings on the general position of CCOP and beyond are as follows:
\begin{itemize}
    \item[(i)] We examine the cardinality-constrained linear independence constraint qualification (CC-LICQ) from \cite{cervinka:2016}. Generically, CC-LICQ is shown to hold at all CCOP feasible points, see Theorem \ref{thm:licq-gen}. Further, we focus on M-stationary points introduced in \cite{burdakov:2016}. For a generic CCOP, M-stationary points turn out to be nondegenerate, see Theorem \ref{thm:generic}.  As in case of nonlinear programming, nondegeneracy refers here to (ND1) CC-LICQ, (ND2) strict complementarity w.r.t active inequality constraints, and (ND4) regularity of Lagrange function's Hessian restricted to a suitably chosen tangent space. As novelty for CCOP, if the cardinality constraint is not active, nondegeneracy additionally requires that (ND3) Lagrange multipliers w.r.t. zero entries of an M-stationary point under consideration do not vanish, see Definition \ref{def:nd-m-st}. In particular, we show that for generic minimizers of CCOP the cardinality constraint has to be always active, see Theorem \ref{thm:gen-min}.
    \item[(ii)] In order to go beyond genericity for CCOP, we study the property of strong stability for M-stationary points. Following the ideas in \cite{kojima:1980}, a strongly stable M-stationary point remains locally unique with respect to any $C^2$-perturbations of the defining functions in $\mbox{CCOP}(f,h,g)$.
    We characterize strongly stable M-stationary points in terms of first- and second-order information of CCOP defining functions, see Theorem \ref{thm:char-ss}. Based on the later,
    nondegeneracy implies strong stability, but not vice versa.
    We emphasize that there could exist strongly stable, but degenerate M-stationary points, -- these are particular singularities of CCOP. Nevertheless, the novel condition ND3 has to be likewise fulfilled at any strongly stable M-stationary point.
    \item[(iii)] Aiming to illustrate the strength of the generic approach, we deal with the global structure of CCOP. For that, we study the topological changes of its lower level sets as their levels vary. Deformation and cell-attachment in the sense of Morse theory, see e.\,g. \cite{milnor:1963}, \cite{goresky:1988}, are proved for generic CCOP. Deformation says that lower level sets are homeomorphic if passing a level which does not correspond to any M-stationary point, see Theorem \ref{thm:def}. Cell-attachment algebraically describes topological differences between lower level sets if a level corresponding to a nondegenerate M-stationary is crossed. It turns out that multiple cells of the same dimension need to be attached to a lower level set in order to obtain another lower level set up to homotopy-equivalence, see Theorem \ref{thm:cell-a}. The dimension of those cells to be attached coincides with the M-index we propose for nondegenerate M-stationary points, see Definition \ref{def:mindex}.  
    A global interpretation of deformation and cell-attachment is given in form of a mounting pass result, see Remark \ref{rem:m-pass}.
\end{itemize}

% Literature CCOP
Let us comment on the relation of our results to those known from the literature on CCOP. The class of CCOP has been introduced in \cite{burdakov:2016}. There, it has been suggested to consider M-stationary points along with the cardinality-constrained linear independence constraint qualification. Subsequent studies in this direction were conducted in \cite{cervinka:2016}, where weaker constraint qualification were examined. In \cite{bucher:2018}, the authors introduce a cardinality-constrained second-order sufficient condition in the framework of CCOP. They show that the fulfilment of this condition at an M-stationary point of CCOP implies its local uniqueness, see Proposition \ref{prop:loc-uniq}. Moreover, if the cardinality constraint is additionally active, the M-stationary point under consideration becomes a strict local minimizer of CCOP, see Proposition \ref{prop:opt-cond}. Another line of research concerning the class of CCOP has been developed in \cite{pan:2017}. There, M-stationary points appear under the name of C-KKT points. By using Bouligand tangential cone w.r.t. the cardinality constraint, an alternative cardinality-constrained second-order sufficient condition has been proposed in \cite{pan:2017}. Under the latter, M-stationary points are shown to be strict local minimizers even if the cardinality constraint becomes inactive, see Proposition \ref{prop:opt-cond1}. However, we emphasize that cardinality-constrained second-order sufficient condition either from \cite{bucher:2018} or from \cite{pan:2017} does not in general guarantee the strong stability of local minimizers, see Examples \ref{ex:inst} and \ref{ex:so-ss}. Throughout the paper, we elaborate in detail on the relation of cardinality-constrained second-order sufficient conditions from \cite{bucher:2018} and from \cite{pan:2017}, respectively, to the notions of nondegenerate, as well as of strongly stable M-stationary points. From the generic viewpoint, in \cite{bucher:2018} and in \cite{pan:2017} some unstable singularities are included into considerations, whereas we concentrate just on nondegeneracy and strong stability for CCOP. 

The paper is organized as follows. In Sections \ref{sec:licq} and \ref{sec:mstat} we discuss the genericity of CC-LICQ, as well as of nondegeneracy of M-stationary points, respectively. Section \ref{sec:ss-m} is devoted to the characterization of strongly stable M-stationary points. In Section \ref{sec:glob}, the global structure of CCOP is described within the scope of Morse theory.

Our notation is standard. The cardinality of a finite set $A$ is denoted by $|A|$. The $n$-dimensional
Euclidean space is denoted by $\R^n$ with the coordinate vectors $e_i,i= 1,\ldots, n$. Its positive orthant is denoted by $\H^n$. Given a twice continuously differentiable function $f:\R^n\rightarrow \R,\nabla f$ denotes its gradient, and $D^2f$ stands for its Hessian. The entries of the subvector $x_I$ correspond to those of $x \in \R^n$ with respect to a given index set $I \subset \{1,\ldots,n\}$.

%%%%%%%%%%%%%%%%%%%%%%%%%%%%%%%%%%%%%%%%%%%%%%%
%Chapter M-Stationarity
%%%%%%%%%%%%%%%%%%%%%%%%%%%%%%%%%%%%%%%%%%%%%%%
\section{Linear independence constraint qualification}
\label{sec:licq}

We shall use the following notation:
\[
\R^{n,k}=\left\{x \in \R^n\, \left\vert\, \left\|x\right\|_0\le k\right.\right\}.
\]
%Occasionally, we shall use the stratification:
%\begin{equation}
%    \label{eq:strat}
%    \R^{n,k}=\bigcup\limits_{\substack{ I \subset \{1,\ldots,n\}\\ |I|\le k}} \bigcup\limits_{J \subset I} Z_{I,J}^n,
%\end{equation}
%where
%\begin{equation}
%    \label{eq:strat2}
%    Z_{I,J}^n=\left\{ x \in \R^n\; \left\vert \; x_{I^c}=0, x_J>0, x_{I\backslash J}<0 \right. \right\}.
%\end{equation}
Moreover, we define the index sets of zero components and of active constraints, respectively:
\[
I_0(x)=\left\{i \in \{1,\ldots, n\}\; \vert\; x_i= 0\right\},\quad 
Q_0(x)=\left\{q \in Q \; \vert\; g_q(x)=0\right\}.
%&I_{g^c}(\bar x)=\left\{i \in \{1,\ldots,m\}\; \vert\; g_i(\bar x)>0\right\}.
\]
 Without loss of generality, we may assume that at a point of interest $x \in M$ with $\left\| x\right\|_0= k$ it holds:
%and $\left|I_{g^c}(\bar x)\right|=\bar m$:
 \[
%%      I_{g^c}(\bar x)&=&\left\{1,\ldots,\bar m\right\},\\
%%      I_{g}(\bar x)&=&\left\{\bar m+1,\ldots,m\right\},\\
       I_{0}(x)=\left\{1,\ldots,k\right\},\quad
      Q_0(x)=\left\{1, \ldots, \left|Q_0(x)\right|\right\}.
\]
We denote $m=|P|+|Q_0(x)|$ and $\ell=|P|+|Q_0(x)|+|I_0(x)|$.

Let us introduce a suitable constraint qualification for CCOP.

\begin{definition}[CC-LICQ, \cite{cervinka:2016}]
\label{def:CC-LICQ}
We say that a feasible point $\bar x \in M$ of CCOP satisfies the cardinality-constrained linear independence constraint qualification (CC-LICQ)  if the following gradients are linearly independent:
\[
\nabla h_p(\bar x), p \in P, \quad
\nabla g_q(\bar x), q \in Q_0(\bar x), \quad
e_i, i \in I_0(\bar x).
\]
\end{definition}

It turns out that CC-LICQ holds generically in the context of CCOP.

\begin{theorem}[Genericity for CC-LICQ]
\label{thm:licq-gen}
Let $\mathcal{F} \subset C^2(\R^n,\R) \times
C^2(\R^n,\R^{|P|}) \times C^2(\R^n,\R^{|Q|})$ be the subset of CCOP defining functions for which all feasible points satisfy CC-LICQ. Then, 
$\mathcal{F}$ is $C_s^2$-open and -dense.
\end{theorem}
\proof
For $x \in \R^n$ we define the so-called $0$-jet:
\[
j^0F(x)=\left(x, h(x), g(x)\right)
\]
of the mapping
\[
  F(x) = \left(h(x), g(x)\right).
\]
Let us define the subset of the image space corresponding to the feasible set of CCOP:
\[
    A=\R^{n,s}\times \{0\}^{|P|} \times \H^{|Q|}.
\]
We use the stratification of
\[
    \label{eq:strat}
    \R^{n,s}=\bigcup\limits_{\substack{ I_1 \subset \{1,\ldots,n\}\\ \left|I_1\right|\le s}} \bigcup\limits_{J_1 \subset I_1} Z_{I_1,J_1},
\]
where
\[
    Z_{I_1,J_1}=\left\{ x \in \R^n\; \left\vert \;  x_{\{1,\ldots,n\} \backslash I_1}=0, x_{J_1}>0, x_{I_1\backslash J_1}<0 \right. \right\}.
\]
We analogously use the stratification of
\[
    \H^{|Q|}=\bigcup\limits_{Q_1 \subset Q} Z_{Q_1},
\]
where
\[
    Z_{Q_1}=\left\{ \left. x \in \R^{|Q|}\; \right\vert \; x_{Q \backslash Q_1}=0, x_{Q_1}>0 \right\}.
\]
The set $A$ admits the following stratification inherited from above:
\[
A=\bigcup\limits_{Q_1 \subset Q}\bigcup\limits_{\substack{ I_1 \subset \{1,\ldots,n\}\\ \left|I_1\right|\le s}} \bigcup\limits_{J_1 \subset I_1} X_{Q_1,I_1,J_1},
\]
where
\[
    X_{Q_1,I_1,J_1}=Z_{I_1,J_1}^n \times \{0\}^p \times Z_{Q_1}.
\]
We show that CC-LICQ equivalently means that $j^0F$ meets $A$ transversally.
Recall that $j^0F$ meets the Whitney stratified set $A$ transversally if for all $\bar x \in j^0F^{-1}(A)$ it holds:
\[
  D j^0 F(\bar x)\left[\R^n\right]+\mathcal{T}_{\bar x}X_{Q_1,I_1,J_1}=\R^{n+|P|+|Q|},
\]
where $X_{Q_1,I_1,J_1}$ is the stratum of $A$ containing $j^0F(\bar x)$.
%Let $\bar x \in j^0F^{-1}(A)$, i.\,e. it is a feasible point of CCOP, and $X_{\bar Q,I,J}$ be the stratum containing $j^0F(\bar x)$.
The differential $Dj^0F(\bar x)\left[\R^n\right]$ is spanned by the columns of the matrix
\[
\left(
\begin{array}{l}
e^T_i, i =1, \ldots, n\\ 
D h_p(\bar x), p \in P\\
D g_q(\bar x), q \in Q\\
\end{array}\right).
\]
The tangent space $\mathcal{T}_{\bar x}X_{Q_1,I_1,J_1}$ is spanned by those unit vectors
from $\R^{n+|P|+|Q|}$, which correspond to the sets $I_1$ and $Q_1$, i.\,e.
\[
\left\{ \left(\begin{array}{c}
     e_i \\ 0 \\ 0
\end{array}\right) \in \R^{n+|P|+|Q|}\;\left\vert\;i\in I_1\right.\right\} \cup \left\{\left(\begin{array}{c}
     0 \\ 0 \\ e_q
\end{array}\right) \in \R^{n+|P|+|Q|}\;\left\vert\; q \in Q_1\right.\right\}.
\]
In order to show that $Dj^0F(\bar x)\left[\R^n\right]$ and $\mathcal{T}_{\bar x}X_{Q_1,I_1,J_1}$ sum up to the whole $\R^{n+|P|+|Q|}$, we determine the column rank of the matrix \[
\left(
\begin{array}{lll}
e^T_i, i =1, \ldots, n&e_i, i \in I_1&0\\ 
D h_p(\bar x), p \in P&0&0\\
D g_q(\bar x), q \in Q&0&e_q, q \in Q_1\\
\end{array}\right).
\]
The column rank is equal to $n+|P|+|Q|$ if and only if all its rows are linearly independent. This is exactly the case if just the following of them -- written as vectors -- are linearly independent:
\[
\left(
\begin{array}{c}
e_i\\ 
0 \\
0
\end{array}\right), i \in \{1,\ldots,n\} \backslash I_1, \quad \left(
\begin{array}{c}
\nabla h_p(\bar x)\\ 
0\\
0
\end{array}\right), p \in P, \quad 
\left(
\begin{array}{c}
\nabla g_q(\bar x) \\ 
0\\
0
\end{array}\right), q \in Q \backslash Q_1.
\]
Hence, $D j^0 F(\bar x)\left[\R^n\right]$ and $\mathcal{T}_{\bar x}X_{Q_1,I_1,J}$ sum up to the whole space $\R^{n+|P|+|Q|}$ if and only if the gradients
\[
e_i, i \in \{1,\ldots,n\} \backslash I_1, \quad \nabla h_p(\bar x), p \in P, \quad
\nabla g_q(\bar x), q \in Q \backslash Q_1
\]
are linearly independent. Since $Q_0(\bar x)=Q \backslash Q_1$ and $I_0(\bar x)=\{1,\ldots,n\} \backslash I_1$, CC-LICQ is valid. 
We apply the structured jet transversality theorem from \cite{guenzel:2008} to conclude the proof. Indeed, the latter says that for a given reduced jet extension and a given stratification the subset of functions, which meet the stratification transversally, is $C^2_s$-dense. For closed stratified sets it also gives that the mentioned subset of functions is $C^1_s$-open. Obviously, it is then also $C_s^2$-open. 
\qed

%\begin{definition}[Coordinate system]
%\label{def:cs}
%The feasible set $M$ admits a local $C^2$-coordinate system of $\R^n$ at $\bar x \in M$ with $\left\|\bar x\right\|_0=k$ by means of a $C^2$-diffeomorhism $\Phi:U \rightarrow V$ with open $\R^n$-neighborhoods $U$ and $V$ of $\bar x$ and $0$, respectively, if it holds that
%\begin{itemize}
%    \item[(i)] $\Phi(\bar x)=0$,
%    \item[(ii)] $\Phi\left(M \cap U\right) = \left(\{0\}^{|P|} \times \H^{\left|Q_0(\bar x)\right|} \times \R^{n-k,s-k} \times \R^{k - |P|-\left|Q_0(\bar x)\right|}\right) \cap V$.
%\end{itemize}
%\end{definition}

Under CC-LICQ, the CCOP feasible set can be locally represented in a product structure by introducing new coordinates.  

\begin{lemma}[Local structure]
\label{lem:local}
Suppose that CC-LICQ holds at $\bar x \in M$ with $\left\|\bar x\right\|_0=k$. Then $M$ admits a local $C^2$-coordinate system of $\R^n$ at $\bar x$, i.\,e. there exists a $C^2$-diffeomorhism $\Phi:U \rightarrow V$ with open $\R^n$-neighborhoods $U$ and $V$ of $\bar x$ and $0$, respectively, such that
\begin{itemize}
    \item[(i)] $\Phi(\bar x)=0$,
    \item[(ii)] $\Phi\left(M \cap U\right) = \left(\{0\}^{|P|} \times \H^{\left|Q_0(\bar x)\right|} \times \R^{n-k,s-k} \times \R^{k - m}\right) \cap V$.
\end{itemize}
\end{lemma}
\proof
Choose vectors $\xi_r\in \R^n, r \in R$,
which together with the vectors 
\[
\nabla h_p(\bar x), p \in P, \quad
\nabla g_q(\bar x), q \in Q_0(\bar x), \quad
e_i, i \in I_0(\bar x)
\]
form a basis for $\R^n$.
We put
\[
\begin{array}{lcl}
    y_p &=&h_p\left( x\right) \mbox{ for } p \in P,  \\
    y_{|P|+q}&=&g_q\left( x\right) \mbox{ for } q \in Q_0(\bar x), \\
    y_{m+i}&=& x_i \mbox{ for } i \in I_0(\bar x),\\
 y_{\ell+r}&=&\xi_r^T\left(x-\bar x\right) \mbox{ for } r \in R.
 \end{array}
\]
We write for short 
\begin{equation}
  \label{eq:stddiff}
y=\Phi(x).  
\end{equation}
By definition it holds $\Phi(x) \in C^2\left(\R^n,\R^n\right)$ and $\Phi(\bar x)=0$. Due to CC-LICQ, the Jacobian matrix $D\Phi(\bar x)$ is nonsingular.
Hence, by means of the inverse function theorem, there exist open neighborhoods
$U$ of $\bar x$ and $V$ of $0$ such that $\Phi:U\rightarrow V$ is a $C^2$-diffeomorphism. Moreover, we can guarantee that $Q_0(x) \subset Q_0(\bar x)$ and $I_0(x) \subset I_0(\bar x)$ by shrinking $U$ if necessary. Thus, property (ii) follows directly from the definition of $\Phi$.\qed

\section{Nondegenerate M-stationary points}
\label{sec:mstat}

Let us introduce a suitable stationarity notion for CCOP.

\begin{definition}[M-stationarity, \cite{burdakov:2016}]
\label{def:M-stat}
A feasible point $\bar x \in M$ is called M-stationary for CCOP if there exist multipliers
$\bar \lambda \in \R^{|P|}$, $\bar \mu \in \R^{\left|Q\right|}$, and $\bar \gamma \in \R^{n}$ such that the following conditions hold:
\begin{itemize}
    \item[] M1: $D f(\bar x) = \sum\limits_{p\in P}\bar \lambda_p D h_p(\bar x)+
    \sum\limits_{q \in Q}\bar \mu_q D g_q(\bar x)
    +\sum\limits_{i =1}^{n} \bar \gamma_i e_i$,
    \item[] M2: $\bar \mu_q g_q(\bar x)= 0$ and $\bar \mu_q \geq 0$ for all $q \in Q$.
   \item[] M3: $\bar \gamma_i \bar x_i=0$ for all $i=1, \ldots, n$.
\end{itemize} 
We call $(\bar x, \bar \lambda, \bar \mu, \bar \gamma)$ an M-stationary pair for CCOP.
\end{definition}
It follows from M2 and M3 that multipliers for non-zero components and inactive constraints vanish:
\[
   \bar \mu_q = 0 \mbox{ for all } q \in Q \backslash Q_0(\bar x), \quad \bar \gamma_i =0 \mbox{ for all } i \in \{1,\ldots, n\} \backslash I_0(\bar x).
\]
Note that under CC-LICQ the multipliers of an M-stationary point are uniquely determined. 
It is straightforward to see that M-stationarity is a first-order sufficient optimality condition for CCOP in this case.

\begin{lemma}[Necessary optimality condition, \cite{burdakov:2016}]
If $\bar x \in M$ is a local minimizer of CCOP satisfying CC-LICQ, then $\bar x$ is M-stationary.
\end{lemma}
Given an $M$-stationary point $\bar x \in M$ with multipliers $(\bar \lambda,\bar \mu,\bar \gamma)$, it is convenient to define the Lagrange function
    \[
     L(x)=f\left(x\right) - \sum\limits_{p \in P}\bar \lambda_p h_p\left( x\right)+
    \sum\limits_{q \in Q}\bar \mu_q g_q\left(x\right)
    +\sum\limits_{i =1}^n \bar \gamma_i x_i.
    \]
Further, we set
\[
  M_0(\bar x) =\left\{
x \in \R^n\,\left\vert\,\begin{array}{l}
      h_p\left(x\right) =0, p \in P \\
g_q\left(x\right)=0,q \in Q_0(\bar x) \\
x_i=0, i \in I_0(\bar x)
\end{array} 
\right.\right\}.
\]
Obviously, $M_0(\bar x) \subset M$. In case that CC-LICQ holds at $\bar x$, the set $M_0(\bar x)$ is locally an  $C^2$-manifold of dimension $\left\|\bar x\right\|_0 - |P|-\left|Q_0(\bar x)\right|$. 
The tangent space of $M_0(\bar x)$ at $\bar x$ is thus given by
\[
    \mathcal{T}_{\bar x} M_0(\bar x)=\left\{
\xi \in \R^n\,\left\vert\, \begin{array}{l} Dh_p(\bar x) \xi=0, p \in P \\
Dg_q(\bar x)\xi=0,q \in Q_0(\bar x)\\
\xi_i=0, i \in I_0(\bar x)
\end{array}
\right.\right\}.
\]

Further, we define the concept of nondegeneracy for M-stationary points.

\begin{definition}[Nondegenerate M-stationarity]
\label{def:nd-m-st}
An M-stationary point $\bar x \in M$ of CCOP is called nondegenerate if the following conditions hold at $\bar x$:
\begin{itemize}
    \item[] ND1: CC-LICQ,
    \item[] ND2: $\bar \mu_q>0$ for all $q\in Q_0(\bar x)$,
    \item[] ND3: if $\left\|\bar x\right\|_0<s$ then $\bar \gamma_i\ne 0$ for all $i\in I_0(\bar x)$,
    \item[] ND4: the matrix $D^2 L(\bar x)\restriction_{\mathcal{T}_{\bar x} M_0(\bar x)}$ is nonsingular.  
    \end{itemize}
\end{definition}

It turns out that nondegeneracy is a generic property of M-stationary points for CCOP.

\begin{theorem}[Genericity for M-stationarity]
\label{thm:generic}
Let $\mathcal{F} \subset C^2(\R^n,\R) \times
C^2(\R^n,\R^{|P|}) \times C^2(\R^n,\R^{|Q|})$ be the subset of CCOP defining functions for which each 
M-stationary point is nondegenerate. Then, $\mathcal{F}$ is $C^2_s$-open and -dense.
\end{theorem}

\proof    Let us fix a number $k\in \left\{0,\ldots,s\right\}$ of non-zero entries, an index set $I_0\subset \left\{1,\ldots,n\right\}$ of $n-k$ zero entries, an index subset $J_0\subset I_0$ of zero-entries, an index set $ Q_0 \subset \left\{1,\ldots,|Q|\right\}$ of active inequality constraints, an index subset $T_0 \subset Q_0$ of these active inequality constraints, and a number $r \in \mathbb{N}$ standing for the rank.
For this choice we consider the set  $M_{k,I_0,J_0,Q_0,T_0,r}$ of $x \in \R^n$ such that the following conditions are satisfied:
\begin{itemize}
    \item[] (m1) $x_i\ne 0$ for all $i \in \{1,\ldots,n\} \backslash I_0$, and $x_i = 0$ for all $i \in I_0$,
    \item[] (m2) $h_p(x)=0$ for all $p \in P$ and $g_{q}(x)=0$ for all $q \in Q_0$,
    \item[] (m3a)
 if $k<s$ then $Df(x) \in \mbox{span}\left\{Dh_p(x), p \in P,
Dg_{q}(x), q \in Q_0\backslash T_0, 
e_i, i \in I_0 \backslash J_0
\right\}$,
    \item[] (m3b) if $k=s$ then $Df(x) \in \mbox{span}\left\{Dh_p(x), p \in P,
D{g_{q}}(x),q \in Q_0\backslash T_0, 
e_i, i \in I_0
\right\}$,
    \item[] (m4) the matrix $D^2 L( x)\restriction_{\mathcal{T}_{x} M_0(x)}$ has rank $r$.
\end{itemize}
Note that (m1) refers to cardinality constraint, (m2) to equality and active inequality constraints, while (m3a) or (m3b) describe violation of ND2 and ND3. Furthermore, (m4) describes violation of ND4.
Now, it suffices to show that $M_{k,I_0,J_0,Q_0,T_0,r}$ is generically empty whenever one of the sets $T_0$ or $J_0$ is nonempty or the rank $r$ in (m4) is not full, i.\,e. $r < \mbox{dim}\left(\mathcal{T}_{x} M_0(x)\right)$.
In fact, the available degrees of freedom
of the variables involved in each $M_{k,I_0,J_0,Q_0,T_0,r}$ are $n$. The loss of freedom caused by (m1) is $n-k$, and
the loss of freedom caused by (m2) is $\left|P\right|+\left|Q_0\right|$.
Due to Theorem \ref{thm:licq-gen}, CC-LICQ holds generically at any feasible $x$, i.\,e. (ND1) is fulfilled. Suppose that the sets $T_0$ and $J_0$ are empty, then both (m3a) or (m3b) causes a loss of freedom of $n-\left|P\right|-\left|Q_0\right|-(n-k)$. Hence, the total loss of freedom is $n$.
We conclude that a
further degeneracy, i.\,e. $T_0 \not = \emptyset$, $J_0 \not = \emptyset$ or $r < \mbox{dim}\left(\mathcal{T}_{x} M_0(x)\right)$, would imply that the total available degrees of freedom $n$ are exceeded. By virtue of the jet
transversality theorem from \cite{Jongen:2000}, generically the sets $M_{k,I_0,J_0,Q_0,T_0,r}$ must be empty.
For the openness result, we argue in a standard way. Locally, M-stationarity can be written
via stable equations. Then, the implicit function theorem for Banach spaces can be applied to
follow M-stationary points with respect to (local) $C^2$-perturbations of defining functions. Finally,
a standard globalization procedure exploiting the specific properties of the strong $C^2_s$-topology can be used to construct a (global) $C^2_s$-neighborhood of problem data for which the nondegeneracy
property is stable, cf. \cite{Jongen:2000}. \qed

As an auxiliary tool, being important for studying the global structure of CCOP later on, we associate with an M-stationary points its M-index.  

\begin{definition}[M-Index]
\label{def:mindex}
Let $\bar x \in M$ with  $\NNorm{\bar x}{0}= k$ be a nondegenerate M-stationary point of CCOP.
The  number of negative eigenvalues of the matrix $D^2 L(\bar x)\restriction_{\mathcal{T}_{\bar x} M_0(\bar x)}$ is called its quadratic index ($QI$). %Furthermore, the number $s-k$ is called its sparsity index ($SI$). 
The number $s-k+QI$ is called the M-index of $\bar x$.
\end{definition}

Let us describe the local structure of CCOP in a sufficiently small neighborhood of a nondegenrate M-stationary point.

\begin{theorem}[Morse Lemma for CCOP]
\label{thm:morse}
Suppose that $\bar x$ is a nondegenerate M-stationary point of CCOP with  $\NNorm{\bar x}{0}= k$ and  quadratic index $QI$. Then, there exist neighborhoods $U_{\bar x}$ and $V_0$ of $\bar x$ and $0$, respectively, and a local $C^1$-coordinate system $\Psi: U_{\bar x} \rightarrow V_0$ of $\R^n$ around $\bar x$ such that:

\begin{equation}
\label{eq:normal}
    f\circ \Psi^{-1}(y)= f(\bar x) +
    \sum\limits_{q \in Q_0(\bar x)} y_{|P|+q} +
    \sum\limits_{i \in I_0(\bar x)}  y_{m+i} + \sum\limits_{r \in R}\pm y_{\ell+r}^2,  
\end{equation}
where $y \in \{0\}^{|P|} \times \H^{\left|Q_0(\bar x)\right|} \times \R^{n-k,s-k} \times \R^{k - m}$. Moreover, there are exactly $QI$ negative squares in (\ref{eq:normal}).
\end{theorem}

\proof Without loss of generality, we may assume $f(\bar x)=0$.
By using $\Phi(x)$ from (\ref{eq:stddiff}), we put $\bar f := f \circ \Phi^{-1}$ on the set 
$\left(\{0\}^{|P|} \times \H^{\left|Q_0(\bar x)\right|} \times \R^{n-k,s-k} \times \R^{k - m} \right)\cap V_0$.
At the origin we have with respect to the new $y$-coordinates:
\begin{itemize}
    \item [(i)] $\displaystyle \frac{\partial \bar f}{\partial y_{|P|+q}} > 0$ for $q \in Q_0(\bar x)$,
    \item [(ii)] if $k<s$ then $\displaystyle \frac{\partial \bar f}{\partial y_{m+i}} \ne 0$ for $i \in I_0(\bar x)$,
    \item [(iii)] $\displaystyle \frac{\partial \bar f}{\partial y_{\ell+r}}=0$ for $r\in R$ and the matrix $\displaystyle \left(\frac{\partial^2 \bar f}{\partial y_{\ell+r_1} \partial y_{\ell+r_2}}\right)_{r_1,r_2 \in R}$ is nonsingular.
\end{itemize}
We denote $\bar f$ by $f$ again. Under the following coordinate transformations the set $$\{0\}^{|P|} \times \H^{\left|Q_0(\bar x)\right|} \times \R^{n-k,s-k} \times \R^{k - m}$$
will be equivariantly transformed in itself. We put $y=\left(Y_{n-r},Y^r\right)$, where 
\[
Y_{n-r}=\left(y_1, \ldots, y_{\ell}\right), \quad
Y^{r}=\left(y_{\ell+1},\ldots,y_n\right).
\]

It holds:
\[
\begin{array}{lcl}
f\left(Y_{n-r},Y^{r}\right)&=&
\displaystyle \int_0^1 \frac{d}{dt}f\left(tY_{n-r},Y^{r}\right)\diff t+f\left(0,Y^r\right)\\ \\
&=&\displaystyle \sum\limits_{q \in Q_0(\bar x)} y_{|P|+q} \cdot d_{|P|+q}(y)+
    \sum\limits_{i \in I_0(\bar x)}  y_{m+i} \cdot d_{m+i}(y) + f\left(0,Y^r\right),
\end{array}
\]
where
\[
\begin{array}{rcl}
   d_{|P|+q}(y)&=&\displaystyle \int_0^1 \frac{\partial f}{\partial y_{|P|+q}}\left(tY_{n-r},Y^{r}\right)\diff t,\quad q \in Q_0(\bar x),
\\ \\
 d_{m+i}(y)&=&\displaystyle \int_0^1 \frac{\partial f}{\partial y_{m+i}}\left(tY_{n-r},Y^{r}\right)\diff t,\quad i \in I_0(\bar x).   
\end{array}
\]
Note that 
 $d_{|P|+q} \in C^1$ for $q \in Q_0(\bar x)$,  and $d_{m+i} \in C^1$ for $i \in I_0(\bar x)$.
Due to (iii), we may apply the standard Morse Lemma on the $C^2$-function $f\left(0,Y^{r}\right)$ without affecting the first $Y_{n-r}$ coordinates, see e.\,g. \cite{Jongen:2000}. The corresponding coordinate transformation is of class $C^1$. Denoting the transformed functions again by $f$ and $d_i$, we obtain
\[
f(y)= \sum\limits_{q \in Q_0(\bar x)} y_{|P|+q} \cdot d_{|P|+q}(y) +
    \sum\limits_{i \in I_0(\bar x)}  y_{m+i} \cdot d_{m+i}(y) + \sum\limits_{r \in R}\pm y_{\ell+r}^2.
\]

In case $k<s$, (i) and (ii) provide that
\[
\begin{array}{rcl}
   d_{|P|+q}(0)&=&\displaystyle \frac{\partial f}{\partial y_{|P|+q}}\left(0\right) > 0,\quad q \in Q_0(\bar x),
\\ \\
 d_{m+i}(0)&=&\displaystyle \frac{\partial f}{\partial y_{m+i}}\left(0\right)\ne 0,\quad i \in I_0(\bar x).   
\end{array}
\]
Hence, we may take
\[
\begin{array}{l}
 y_{|P|+q} \cdot \left|d_{|P|+q}(y)\right|, \quad q \in Q_0(\bar x), \\
  y_{m+i} \cdot d_{m+i}(y), \quad  i \in I_0(\bar x), \\
y_{\ell+r}, \quad r \in R
\end{array}
\]
as new local $C^1$-coordinates by a straightforward application of the inverse function theorem. Denoting the transformed function again by $f$, we obtain (\ref{eq:normal}). Here, the coordinate transformation $\Psi$
is understood as the composition of all previous ones.

In case $k=s$, we need to consider $f$ locally around the origin on the set
\[
\displaystyle\{0\}^{|P|} \times \H^{\left|Q_0(\bar x)\right|} \times \R^{n-k,s-k} \times \R^{k - m}
=
\{0\}^{|P|} \times \H^{\left|Q_0(\bar x)\right|} \times\{0\}^{n-s} \times \R^{k - m}.
\]
Hence, $y_{m+i}=0$ for $i \in I_0(\bar x)$ and we obtain the representation (\ref{eq:normal}) analogously. \qed 

\begin{remark}
\label{rem:multipliers}
It follows from the proof of Lemma \ref{lem:local} and from Lemma 2.2.1 of \cite{Jongen:2004} that the multipliers at a nondegenerate M-stationary point are the corresponding partial derivatives of the objective function in new coordinates given by the diffeomorphism (\ref{eq:stddiff}).
\end{remark}
As a first application of Morse Lemma we analytically describe nondegenerate local minimizers.  

\begin{lemma}[Minimizers and M-index]
\label{lem:min-index}
Let $\bar x \in M$ be a nondegenerate M-stationary point. Then, $\bar x$ is a local
minimizer of CCOP if and only if its M-index vanishes.
\end{lemma}

\proof Let $\bar x$ be a nondegenerate M-stationary point for CCOP.
The application of Morse Lemma from Theorem \ref{thm:morse} says that there exist
neighborhoods $U_{\bar x}$ and $V_0$ of $\bar x$ and $0$, respectively, and a local $C^1$-coordinate system $\Psi: U_{\bar x} \rightarrow V_0$ of $\R^n$ around $\bar x$ such that (\ref{eq:normal}) holds.
Therefore, $\bar x$ is a local minimizer for CCOP if and only if 0 is a local minimizer of $ f\circ \Psi^{-1}$ on the set $\{0\}^{|P|} \times \H^{\left|Q_0(\bar x)\right|} \times \R^{n-k,s-k} \times \R^{k - m}$.
If the M-index vanishes, we have $k=s$ and $QI=0$, and (\ref{eq:normal}) reads as
\begin{equation}
    \label{eq:normalatMI0}
f\circ \Psi^{-1}(y)= f(\bar x) +
    \sum\limits_{q \in Q_0(\bar x)} y_{|P|+q} +
    \sum\limits_{r \in R} y_{\ell+r}^2,
    \end{equation}
where $y \in \{0\}^{|P|} \times \H^{\left|Q_0(\bar x)\right|} \times \{0\}^{n-s} \times \R^{k - m}$. Thus, $0$ is a local minimizer for (\ref{eq:normalatMI0}).
Vice versa, if $0$ is a local minimizer for (\ref{eq:normal}), then obviously $k=s$ and $QI= 0$, hence, the M-index of $\bar x$ vanishes.
\qed

Lemma \ref{lem:min-index} motivates to introduce the notion of nondegenerate local minimizers.

\begin{definition}[Nondegenerate minimizers]
\label{def:nd-min}
A local minimizer $\bar x \in M$ of CCOP is called nondegenerate if the following conditions hold at $\bar x$:
\begin{itemize}
    \item[] CC-LICQ,
    \item[] Strict Complementarity (SC), i.\,e. $\bar \mu_q>0$ for all $q\in Q_0(\bar x)$,
    \item[] Active Cardinality Condition (ACC), i.\,e. $\left\|\bar x\right\|_0=s$,
    \item[] Second-Order Sufficiency Condition (SOSC), i.\,e. $D^2 L(\bar x)\restriction_{\mathcal{T}_{\bar x} M_0(\bar x)}$ is positive definite.  
    \end{itemize}
\end{definition}

Let us relate the notion of nondegenerate minimizer to the cardinality-constrained second-order sufficient conditions from \cite{bucher:2018}. In order to formulate the latter, a substitute for the linearization cone of $M$ at a feasible point $\bar x \in M$ is used:
\[
\mathcal{L}_{\bar x} M = \left\{ \xi \in \R^n \, \left|\, \begin{array}{l}
      Dh_p(\bar x)\xi =0, p \in P \\ Dg_q(\bar x) \xi \geq 0,q \in Q_0(\bar x) \\ \left|\left\{ i \in I_0(\bar x) \,\left|\, \xi_i =0\right.\right\}\right| \geq n-s
\end{array}\right.\right\}.
\]
The critical cone of $M$ at $\bar x \in M$ is then defined by
\[
  \mathcal{C}_{\bar x} M = \left\{ \xi \in \mathcal{L}_{\bar x} M \, \left|\, Df(\bar x) \xi \leq 0 \right.\right\}.
\]
Next Lemma \ref{lem:crit-cone} provides a representation of the critical cone just in terms of the constraints at an M-stationary point under ACC.

\begin{lemma}[Critical cone] 
\label{lem:crit-cone}
Let $\bar x \in M$ be an M-stationary point of CCOP with multipliers $(\bar \lambda,\bar \mu,\bar \gamma)$ satisfying ACC. Then, for the critical cone of $M$ at $\bar x$ holds: 
\begin{equation}
    \label{eq:cc-rep}
    \mathcal{C}_{\bar x} M = \left\{ \xi \in \R^n \, \left|\, 
    \begin{array}{l}
          Dh_p(\bar x)\xi =0, p \in P \\ Dg_q(\bar x) \xi \geq 0,q \in Q_0(\bar x) \backslash  Q_+(\bar x), Dg_q(\bar x) \xi = 0, q \in Q_+(\bar x) \\
          \xi_i=0, i \in I_0(\bar x)
    \end{array}
\right.\right\},
\end{equation}
where the index set of positive multipliers corresponding to the inequality constraints is given by
\[
Q_+(\bar x) = \left\{ q \in Q_0(\bar x) \,\left|\, \bar \mu_q > 0\right.\right\}.
\]
\end{lemma}

\proof  Let us consider a vector $\xi$ from the right-hand side of (\ref{eq:cc-rep}). Due to ACC, we have $\left|I_0(\bar x) \right|=n-s$. Hence, from $\xi_i=0, i \in I_0(\bar x)$ we get:
\[
\left|\left\{ i \in I_0(\bar x) \,\left|\, \xi_i =0\right.\right\}\right| = n-s.
\]
This implies that $\xi \in \mathcal{L}_{\bar x} M$.
M-stationarity of $\bar x$ provides:
\[
 \begin{array}{lcl}
      D f(\bar x) \xi &=& \displaystyle \sum\limits_{p\in P}\bar \lambda_p \underbrace{D h_p(\bar x) \xi}_{=0} + \sum\limits_{q \in Q_0(\bar x) \backslash Q_+(\bar x)} \underbrace{\bar \mu_q}_{=0} D g_q(\bar x) \xi \\ \\ &&\displaystyle +
    \sum\limits_{q \in Q_+(\bar x)}\bar \mu_q \underbrace{D g_q(\bar x) \xi}_{=0}
    +\sum\limits_{i \in I_0(\bar x)} \bar \gamma_i \underbrace{e_i \xi}_{= \xi_i=0} = 0.  
 \end{array}
\]
Overall, we have shown that $\xi \in \mathcal{C}_{\bar x} M$. Now, we assume that $\xi \in \mathcal{C}_{\bar x} M$. By recalling $\left|I_0(\bar x) \right| = n-s$, the condition  $\left|\left\{ i \in I_0(\bar x) \,\left|\, \xi_i =0\right.\right\}\right| \geq n-s$ implies:
\[
\xi_i=0, i \in I_0(\bar x).
\]
Furthermore, we obtain: 
\[
 \begin{array}{lcl}
   \displaystyle \sum\limits_{q \in Q_+(\bar x)}\underbrace{\bar \mu_q}_{>0} \underbrace{D g_q(\bar x) \xi}_{\geq 0} &=& \displaystyle \underbrace{D f(\bar x) \xi}_{\leq 0} -  \sum\limits_{p\in P}\bar \lambda_p \underbrace{D h_p(\bar x) \xi}_{=0} \\ \\ && \displaystyle -\sum\limits_{q \in Q_0(\bar x) \backslash Q_+(\bar x)} \underbrace{\bar \mu_q}_{=0} D g_q(\bar x) \xi -
    \sum\limits_{i \in I_0(\bar x)} \bar \gamma_i \underbrace{e_i \xi}_{= \xi_i=0}.
 \end{array}
\]
From here we deduce that $D g_q(\bar x) \xi=0$ for all $q \in Q_0(\bar x)$. \qed

The following related sufficient optimality condition has been stated as Corollary 3.2 in \cite{bucher:2018}.

\begin{proposition}[Sufficient optimality condition, \cite{bucher:2018}]
\label{prop:opt-cond}
Let  $\bar x \in M$ be an M-stationary point satisfying ACC. Assume that for all $\xi \in \mathcal{C}_{\bar x} M$ with $\xi \not =0$ there exist multipliers $\bar \lambda \in \R^{|P|}$ and $\bar \mu \in \R^{\left|Q\right|}$ such that:
\begin{equation}
    \label{eq:cc-sosc32}
    \xi^T \left(D^2f\left(\bar x\right) - \sum\limits_{p \in P}\bar \lambda_p D^2h_p\left(\bar x\right)-
    \sum\limits_{q \in Q}\bar \mu_q D^2g_q\left(\bar x\right)\right) \xi >0.
\end{equation}
Then, $\bar x$ is a strict local minimizer of CCOP, i.\,e. there exists $r > 0$ such that $f(\bar x)<f(x)$ holds for all $x \in M \cap B\left(\bar x, r \right)$ with $x \not= \bar x$.
\end{proposition}

It is straightforward to apply Proposition \ref{prop:opt-cond} for M-stationary points $\bar x \in M$, which satisfy CC-LICQ, SC, ACC, and SOSC. In this case, CC-LICQ ensures the uniqueness of Lagrange multipliers $(\bar \lambda, \bar \mu, \bar \gamma)$. SC implies that $Q_+(\bar x) = Q_0(\bar x)$. Due to ACC, Lemma \ref{lem:crit-cone} is applicable and we obtain $\mathcal{C}_{\bar x} M = \mathcal{T}_{\bar x} M_0(\bar x)$. Hence, the cardinality-constrained second-order sufficient condition (\ref{eq:cc-sosc32}) from \cite{bucher:2018} coincides with SOSC.
Altogether, the assumptions of Proposition \ref{prop:opt-cond} are fulfilled and it follows that $\bar x$ is a strict local minimizer. We conclude that the notion of
nondegenerate minimizer is in accordance with the cardinality-constrained second-order sufficient conditions from \cite{bucher:2018}. More precisely, Lemma \ref{lem:min-index} can be partly deduced by means of Proposition \ref{prop:opt-cond}, i.\,e. that a nondegenerate M-stationary point with vanishing M-index (thus, satisfying CC-LICQ, SC, ACC, and SOSC) is a local minimizer of CCOP.
From this point of view, our contribution here is not so much in proving Lemma \ref{lem:min-index}, but rather in recognizing that all local minimizers are generically nondegenerate.

\begin{theorem}[Genericity for minimizers]
\label{thm:gen-min}
Let $\mathcal{F} \subset C^2(\R^n,\R) \times
C^2(\R^n,\R^{|P|}) \times C^2(\R^n,\R^{|Q|})$ be the subset of CCOP defining functions for which each 
local minimizer is nondegenerate, i.\,e. satisfying CC-LICQ, SC, ACC, and SOSC. Then, $\mathcal{F}$ is $C^2_s$-open and -dense.
\end{theorem}

\proof Note that every local minimizer of CCOP has to be M-stationary.
Nondegenarate M-stationary points are generic by Theorem \ref{thm:generic}. Hence, generically, local minimizers are nondegenarate. Thus, CC-LICQ and SC are satisfied. Note that moreover its M-index vanishes due to Lemma \ref{lem:min-index}. Hence, ACC and SOSC are also satisfied. \qed

Additionally, we would like to relate the notion of nondegenerate minimizer to the cardinality-constrained second-order sufficient conditions from \cite{pan:2017}. For that, let $\bar x \in M$ be an M-stationary point for CCOP with multipliers $(\bar \lambda,\bar \mu,\bar \gamma)$. A corresponding linearization cone w.r.t. equality and inequality constraints at $\bar x \in M$ is defined in \cite{pan:2017} as follows:
\begin{equation}
    \label{eq:cc-rep1}
    \mathcal{L}_{Q_+}(\bar x) = \left\{ \xi \in \R^n \, \left|\, 
    \begin{array}{l}
          Dh_p(\bar x)\xi =0, p \in P \\ Dg_q(\bar x) \xi \geq 0,q \in Q_0(\bar x) \backslash  Q_+(\bar x), Dg_q(\bar x) \xi = 0, q \in Q_+(\bar x)
    \end{array}
\right.\right\}.
\end{equation}
Moreover, the authors compute the Bouligand tangential cone w.r.t. the cardinality constraint:
\begin{equation}
    \label{eq:cc-b}
  T^B_{\R^{n,s}}(\bar x)= \left\{
  \begin{array}{cc}
     \displaystyle \bigcup_{J \in \mathcal{J}(\bar x)} \mbox{span}\left\{ e_j \,\left|\, j \in J\right.\right\}, & \mbox{if } \|\bar x\|_0<s, \\
  \mbox{span}\left\{ e_i \,\left|\, i \in I_1(\bar x)\right.\right\}, & \mbox{if } \|\bar x\|_0=s, \\
  \end{array}
\right.
\end{equation}
where
\[
  \mathcal{J}(\bar x) =\left\{ J \subset \{1,\ldots,n\} \,\left|\, I_1(\bar x) \subset J, |J|=s\right.\right\}.
\]

The following related sufficient optimality condition has been stated as Theorem 4.2 in \cite{pan:2017}.

\begin{proposition}[Sufficient optimality condition, \cite{pan:2017}]
\label{prop:opt-cond1}
Let  $\bar x \in M$ be an M-stationary point with multipliers $\bar \lambda \in \R^{|P|}$ and $\bar \mu \in \R^{\left|Q\right|}$. Assume that for all $\xi \in \mathcal{L}_{Q_+}(\bar x) \cap T^B_{\R^{n,s}}(\bar x)$ with $\xi \not =0$ it holds:
\begin{equation}
    \label{eq:cc-sosc321}
    \xi^T \left(D^2f\left(\bar x\right) - \sum\limits_{p \in P}\bar \lambda_p D^2h_p\left(\bar x\right)-
    \sum\limits_{q \in Q}\bar \mu_q D^2g_q\left(\bar x\right)\right) \xi >0.
\end{equation}
Then, $\bar x$ is a strict local minimizer of CCOP, i.\,e. there exists $r > 0$ such that $f(\bar x)<f(x)$ holds for all $x \in M \cap B\left(\bar x, r \right)$ with $x \not= \bar x$.
\end{proposition}

In view of Lemma \ref{lem:crit-cone}, it is now straightforward to see that under ACC at an M-stationary point $\bar x \in M$ with multipliers $(\bar \lambda,\bar \mu,\bar \gamma)$ we have:
\begin{equation}
    \label{eq:under-acc}
   \mathcal{C}_{\bar x} M = \mathcal{L}_{Q_+}(\bar x) \cap T^B_{\R^{n,s}}(\bar x).
\end{equation}
Hence, cardinality-constrained sufficient optimality conditions (\ref{eq:cc-sosc32}) from \cite{bucher:2018} and (\ref{eq:cc-sosc321}) from \cite{pan:2017} coincide if ACC additionally holds. From here we conclude that the notion of
nondegenerate minimizer is  in accordance also with the cardinality-constrained second-order sufficient conditions from \cite{pan:2017}. Again, Lemma \ref{lem:min-index} can be partly deduced by means of Proposition \ref{prop:opt-cond1}, i.\,e. that a nondegenerate M-stationary point with vanishing M-index (thus, satisfying CC-LICQ, SC, ACC, and SOSC) is a local minimizer of CCOP. From the generic viewpoint, the authors in \cite{bucher:2018} and in \cite{pan:2017} consider not only nondegenerate M-stationary points, but also some singularities which satisfy
less demanding cardinality-constrained second-order sufficient conditions (\ref{eq:cc-sosc32}) and (\ref{eq:cc-sosc321}), respectively. However, as we shall see in next Section \ref{sec:ss-m}, degenerate M-stationary points may become unstable w.r.t. data perturbations, even in presence of the cardinality-constrained second-order sufficient condition (\ref{eq:cc-sosc32}) or (\ref{eq:cc-sosc321}).

\section{Strongly stable M-stationary points}
\label{sec:ss-m}

In Corollary 3.3 by \cite{bucher:2018}, the local uniqueness of M-stationary points has been deduced in terms of cardinality-constrained second-order sufficient conditions. In order to state the corresponding result, we mention the constant positive linear dependence constraint qualification.

\begin{definition}[CC-CPLD, \cite{cervinka:2016}]
\label{def:CC-CPLD}
We say that a feasible point $\bar x \in M$ of CCOP satisfies the cardinality-constrained constant positive linear dependence constraint qualification (CC-CPLD) if for any subset $\bar Q_0 \subset Q_0(\bar x)$, $\bar P \subset P$, and $\bar I_0 \subset I_0(\bar x)$, such that the gradients
\[
\nabla h_p(x), p \in \bar P, \quad
\nabla g_q( x), q \in \bar Q_0, \quad
e_i, i \in \bar I_0
\]
are positively linearly dependent at $x=\bar x$, they remain linearly dependent in a neighborhood of $\bar x$.
\end{definition}
It is not hard to see from Definitions \ref{def:CC-LICQ} and \ref{def:CC-CPLD} that CC-LICQ implies CC-CPLD, see \cite{cervinka:2016}.

\begin{proposition}[Local uniqueness, \cite{bucher:2018}]
\label{prop:loc-uniq}
Let  $\bar x \in M$ be an M-stationary point of CCOP satisfying CC-CPLD. Assume that for all $\xi \in \mathcal{C}_{\bar x} M$ with $\xi \not =0$ and all multipliers $\bar \lambda \in \R^{|P|}$ and $\bar \mu \in \R^{\left|Q\right|}$ it holds:
\begin{equation}
    \label{eq:cc-sosc33}
       \xi^T \left(D^2f\left(\bar x\right) - \sum\limits_{p \in P}\bar \lambda_p D^2h_p\left(\bar x\right)-
    \sum\limits_{q \in Q}\bar \mu_q D^2g_q\left(\bar x\right)\right) \xi >0.
\end{equation}
Then, there exists $r>0$ such that $\bar x$ is the unique M-stationary point  within the ball $B(\bar x, r)$.
\end{proposition}

We intend to consider a more demanding property of M-stationary points, namely that of strong stability in the sense of \cite{kojima:1980}. Loosely speaking, a strongly stable M-stationary point remains locally unique with respect to any $C^2$-perturbations of the defining functions in $\mbox{CCOP}(f,h,g)$. In order to control these $C^2$-perturbations, we use the seminorm $\left\|(f,h,g)\right\|^{C^2}_{B(\bar x, r)}$ to be the modulus of the function values and partial derivatives up to order two of $f,h,g$ on the ball $B(\bar x, r)$.

\begin{definition}[Strongly stable M-stationary point]
\label{def:s-stab}
An M-stationary point $\bar x$ of $\mbox{CCOP}(f,h,g)$ is called strongly stable if for some $r>0$ and each $\varepsilon \in (0,r]$ there exists $\delta >0$ such that whenever
\[
  \left(\widetilde f, \widetilde h, \widetilde g \right) \in C^2(\R^n,\R) \times
C^2(\R^n,\R^{|P|}) \times C^2(\R^n,\R^{|Q|})
\]
and
\[
\left\|\left(\widetilde f, \widetilde h, \widetilde g \right) - (f,h,g)\right\|^{C^2}_{B(\bar x, r)} \leq \delta,
\]
the ball $B\left(\bar x, \varepsilon \right)$ contains an M-stationary point $\widetilde x$ of $\mbox{CCOP}\left(\widetilde f, \widetilde h, \widetilde g \right)$ that is unique within $B\left(\bar x, r \right)$.
\end{definition}

It turns out that the cardinality-constrained second-order sufficient condition (\ref{eq:cc-sosc33}) from \cite{bucher:2018} cannot in general prevent an M-stationary point from being unstable. By the way, the same is true for the cardinality-constrained second-order sufficient condition (\ref{eq:cc-sosc321}) from \cite{pan:2017}.
 
\begin{example}[Instability]
\label{ex:inst}
We consider the following CCOP with $P=Q=\emptyset$ and $n=2$, $s=1$:
\[
 \mbox{CCOP}(f): \quad \min_{x_1,x_2}\,\, f\left(x_1,x_2\right)= x_1^2 + x_2^2 \quad \mbox{s.\,t.} \quad 
   \left\|\left(x_1, x_2\right)\right\|_0 \leq 1.
\]
Obviously, $\bar x=(0,0)$ is the unique minimizer of $\mbox{CCOP}(f)$. Since CC-LICQ is satisfied at $\bar x$, so is also CC-CPLD. Moreover, we have:
\[
\mathcal{C}_{\bar x} M=\left\{(\xi_1,\xi_2) \in \R^2 \,\left|\, \left\|\left(\xi_1, \xi_2\right)\right\|_0 \leq 1\right. \right\}.
\]
The cardinality-constrained second-order sufficient condition (\ref{eq:cc-sosc33}) from \cite{bucher:2018} is valid, since for all $\xi \in \mathcal{C}_{\bar x} M$ with $\xi \not = (0,0)$ it holds:
\[
  \xi^T D^2f(\bar x) \xi = 2 \left(\xi_1^2 + \xi_2^2\right)>0.
\]
Although both assumptions from Proposition \ref{prop:loc-uniq} are satisfied, $\bar x$ is not strongly stable. To see this, we perturb the defining function by means of an arbitrarily small $\varepsilon > 0$ as follows:
\[
\mbox{CCOP}\left(\widetilde f\right):\quad \min_{x_1,x_2}\,\, \widetilde f\left(x_1,x_2\right)=\left(x_1-\varepsilon\right)^2 + \left(x_2-\varepsilon\right)^2 \quad \mbox{s.\,t.} \quad 
   \left\|\left(x_1, x_2\right)\right\|_0 \leq 1.
\]
Obviously, $\mbox{CCOP}\left(\widetilde f\right)$ has now two solutions $\widetilde x^{1}=(\varepsilon,0)$ and $\widetilde x^{2}=(0,\varepsilon)$. 
Here, we observe a bifurcation of the minimum $\bar x$ of the original problem $\mbox{CCOP}(f)$ into two minima $\widetilde x^{1}$ and $\widetilde x^{2}$ of the perturbed problem $\mbox{CCOP}\left(\widetilde f\right)$.  More interestingly, there is another M-stationary point $\widetilde x^3 =(0,0)$ of the perturbed problem in an arbitrarily small neighborhood of $\bar x$. 
Why does the bifurcation occur? This is not only due to the fact the minimizer $\bar x$ of $\mbox{CCOP}(f)$ is degenerate, but mainly because ACC is violated at $\bar x$, i.\,e. $\left\| \bar x\right\|_0=0$. It is worth to mention that the cardinality-constrained second-order sufficient condition (\ref{eq:cc-sosc321}) from \cite{pan:2017} is nevertheless valid at $\bar x$. In absence of equality and inequality constraints, we have $\mathcal{L}_{Q_+}(\bar x) =\R^2$. It is easy to see that additionally $T^B_{\R^{2,1}}(\bar x) = C_{\bar x} M$ holds here. We conclude that the cardinality-constrained second-order sufficient condition (\ref{eq:cc-sosc321}) from \cite{pan:2017}, although fulfilled at $\bar x$, does not prevent the latter M-stationary point from being unstable.
\qed
\end{example}

Now, we are ready to fully characterize the strong stability in the context of CCOP under CC-LICQ. For that, we shall use some auxiliary objects associated with an M-stationary point $\bar x \in M$ and its mulitpliers $(\bar \lambda,\bar \mu, \bar \gamma)$.
For $Q_+(\bar x) \subset Q_*\subset Q_0(\bar x)$ we set 
\[
  M_*(\bar x) =\left\{
x \in \R^n\,\left\vert\, \begin{array}{l}
    h_p\left(x\right) =0, p \in P\\
g_q\left(x\right)=0,q \in Q_*\\
x_i=0, i \in I_0(\bar x)  
\end{array}
\right.\right\}.
\]
Obviously, $M_*(\bar x) \subset M$. In case that CC-LICQ holds at $\bar x$, the set $M_*(\bar x)$ is locally an  $C^2$-manifold of dimension $\left\|\bar x\right\|_0 - |P|-\left|Q_*\right|$. 
The tangent space of $M_*(\bar x)$ at $\bar x$ is thus given by
\[
\mathcal{T}_{\bar x} M_*(\bar x)=\left\{
\xi \in \R^n\,\left\vert\, \begin{array}{l}
   Dh_p(\bar x) \xi=0, p \in P\\
Dg_q(\bar x)\xi=0,q \in Q_*\\
\xi_i=0, i \in I_0(\bar x)
\end{array}
\right.\right\}.
\]

\begin{theorem}[Characterization of strong stability]
\label{thm:char-ss}
Let $\bar x \in M$ be an $M$-stationary point of CCOP satisfying CC-LICQ. Then, $\bar x$ is strongly stable if and only if it fulfils ND3 and for all index subsets $Q_+(\bar x) \subset Q_* \subset Q_0(\bar x)$ the matrices $D^2 L(\bar x)\restriction_{\mathcal{T}_{\bar x} M_*(\bar x)}$ are nonsingular with the same determinant sign. 
\end{theorem}

In order to prove Theorem \ref{thm:char-ss}, we need the notion of strongly stable M-stationary pairs and their characterization given below.

\begin{definition}[Strongly stable M-stationary pair]
\label{def:s-stab-pair}
An M-stationary pair $\bar x$ of $\mbox{CCOP}(f,h,g)$ along with the corresponding multipliers $(\bar \lambda, \bar \mu, \bar \gamma)$ is called strongly stable if for some $r>0$ and each $\varepsilon \in (0,r]$ there exists $\delta >0$ such that whenever
\[
  \left(\widetilde f, \widetilde h, \widetilde g \right) \in C^2(\R^n,\R) \times
C^2(\R^n,\R^{|P|}) \times C^2(\R^n,\R^{|Q|})
\]
and
\[
\left\|\left(\widetilde f, \widetilde h, \widetilde g \right) - (f,h,g)\right\|^{C^2}_{B(\bar x, r)} \leq \delta,
\]
the ball $B\left((\bar x, \bar \lambda, \bar \mu, \bar \gamma), \varepsilon \right)$ contains an M-stationary pair $\left(\widetilde x, \widetilde \lambda, \widetilde \mu, \widetilde \gamma\right)$ of $\mbox{CCOP}\left(\widetilde f, \widetilde h, \widetilde g \right)$ that is unique within $B\left((\bar x, \bar \lambda, \bar \mu, \bar \gamma), r \right)$.
\end{definition}

\begin{lemma}
\label{lem:stablepair}
The following assertions are equivalent:
\begin{itemize}
    \item [(i)] $\bar x$ is a strongly stable M-stationary point for CCOP which statisfies CC-LICQ and has the associated multiplier vector
    $\left(\bar \lambda, \bar \mu, \bar \gamma\right)$.
    \item [(ii)]$\left(\bar x, \bar \lambda, \bar \mu, \bar \gamma\right)$ is a strongly stable M-stationary pair for CCOP.
\end{itemize}
\end{lemma}
\proof 
\noindent The proof goes along the lines of \cite{klatte:1990}.

(i) $\Rightarrow$ (ii) Consider a strongly stable M-stationary point $\bar x$ fulfilling CC-LICQ and
having the associated multiplier vector
    $\left(\bar \lambda, \bar \mu, \bar \gamma\right)$.
CC-LICQ ensures the uniqueness of $\left(\bar \lambda, \bar \mu, \bar \gamma\right)$. Moreover, the CC-LICQ remains valid under small perturbations. Thus, uniqueness of the Lagrange-multipliers is provided under small perturbations. Moreover, due to Lemma \ref{lem:local} the Langrange-Multipliers are the corresponding partial derivatives of the objective function in new coordinates.
Hence, continuity under small perturbations is also provided.
We conclude, that $\left(\bar x, \bar \lambda, \bar \mu, \bar \gamma\right)$ is a strongly stable M-stationary pair for CCOP.

(i) $\Rightarrow$ (ii) 
Suppose now $\left(\bar x, \bar \lambda, \bar \mu, \bar \gamma\right)$ is a strongly stable M-stationary pair for $\mbox{CCOP}(f,h,g)$. From the definition of strong stability it follows trivially that $\left(\bar x \right)$ is a strongly stable M-stationary point for CCOP with its Lagrange multiplier vector being $\left(\bar \lambda, \bar \mu, \bar \gamma\right)$. It remains to show that it also satisfies CC-LICQ.
In order to show this we suppose LICQ to be not fulfilled. Thus, there exist numbers $\beta_{h,p},p \in P,\beta_{g,q},q\in Q_0\left(\bar x\right),\beta_{x,i}, i \in I_0\left(\bar x\right)$ (not all vanishing) such that
\begin{equation}
\label{eq:lindependence}
    \sum\limits_{p\in P}\bar\beta_{h,p} D h_p(\bar x)+
    \sum\limits_{q \in Q_0\left(\bar x \right)}\bar \beta_{g,q} D g_q(\bar x)
    +\sum\limits_{i\in I_0\left(\bar x\right)} \bar \beta_{x,i} e_i=0.
\end{equation}
Next, we define
\[
c=-\left( \sum\limits_{p\in P} D h_p(\bar x)+
    \sum\limits_{q \in Q_0\left(\bar x \right)} D g_q(\bar x)
    +\sum\limits_{i\in I_0\left(\bar x\right)}  e_i\right)
\]
and $\varphi(x)=c^T\cdot x$.
For $\epsilon>0$ we put:
\[
\begin{array}{lcl}
     \lambda_p(\epsilon)&=&\bar \lambda_p + \epsilon,\quad p\in P,\\ \\
     \mu_q(\epsilon)&=& \begin{cases}
     0,& q \in Q_0^c\left(\bar x \right)\\ 
     \bar \mu_q + \epsilon, &q \in Q_0^c\left(\bar x \right),\end{cases}\\ \\
     \gamma_i(\epsilon)&=& \begin{cases}
     0,& i \in I_1\left(\bar x\right)\\
     \bar \gamma_i + \epsilon,& i \in I_0\left(\bar x\right).
     \end{cases}
\end{array}
\]
Hence, it obviously holds:
\begin{itemize}
    \item$ \displaystyle D f(\bar x) - \epsilon \cdot c = \sum\limits_{p\in P} \lambda_p(\epsilon) D h_p(\bar x)+
    \sum\limits_{q \in Q}\mu_q(\epsilon) D g_q(\bar x)
    +\sum\limits_{i =1}^{n} \gamma_i(\epsilon) e_i$,
    \item$\displaystyle \mu_q(\epsilon) g_q(\bar x)= 0$ and $\mu_q(\epsilon) \geq 0$ for all $q \in Q$.
   \item $\displaystyle \gamma_i(\epsilon) \bar x_i=0$ for all $i=1, \ldots, n$.
\end{itemize}
Hence, $\left(\bar x, \lambda(\epsilon),\mu(\epsilon),\gamma(\epsilon)\right) $ is
a M-stationary pair for $\mbox{CCOP}(f+\epsilon \varphi,h,g)$. For  sufficiently small $\epsilon$ the M-stationary pair$\left(\bar x, \lambda(\epsilon),\mu(\epsilon),\gamma(\epsilon)\right)$
has to be unique in some neighborhood $U$ of $\left(\bar x, \bar \lambda, \bar \mu, \bar \gamma\right)$. However, (\ref{eq:lindependence}) and $\mu_q(\epsilon) >0, q \in Q_0^c\left(\bar x \right)$, ensure together that there is another M-stationary pair for $\mbox{CCOP}(f+\epsilon \varphi,h,g)$ given by$
\left(\bar x, \tilde\lambda(\epsilon,t, \beta), \tilde\mu(\epsilon,t, \beta) , \tilde\gamma(\epsilon,t, \beta) \right)$, with

\[
\begin{array}{lcl}
     \tilde\lambda(\epsilon,t, \beta)&=&\lambda_p(\epsilon)+t\cdot\beta_{h,p}, \quad p\in P,\\ \\
     \tilde\mu(\epsilon,t, \beta)&=& \begin{cases}
     0,& q \in Q_0^c\left(\bar x \right)\\
    \mu_q(\epsilon)+t \cdot \beta_{g,q},&q \in Q_0\left(\bar x \right),\end{cases}\\ \\
     \tilde\gamma(\epsilon,t, \beta)&=& \begin{cases}
     0,& i \in I_1\left(\bar x\right)\\
   \gamma_i(\epsilon)+ t\cdot \beta_{x,i},& i \in I_0\left(\bar x\right),
     \end{cases}
\end{array}
\]
and $t$ a sufficiently small real number. Thus, $\left(\bar x, \bar \lambda, \bar \mu, \bar \gamma\right)$ is not a strongly stable M-stationary pair for $\mbox{CCOP}(f,h,g)$, which is contradictory. Hence, CC-LICQ has to hold at $\bar x$.

\qed

\noindent
{\bf Proof of Theorem \ref{thm:char-ss}:}

First, we consider the \emph{necessity part} and assume that CC-LICQ holds and the M-stationary point $\bar x$ fulfils ND3. Due to Lemma \ref{lem:stablepair}, we may instead observe the M-stationary pair $\left(\bar x, \bar \lambda, \bar \mu, \bar \gamma \right)$.
Now, for any sufficiently close pair $\left( x, \lambda, \mu, \gamma \right)$  it holds $I_0\left( x\right)\subset I_0\left(\bar x\right)$ due to continuity arguments. Moreover, we claim that the sets will be equal, i.\,e. $I_0\left( x\right)= I_0\left(\bar x\right)$, if $\left( x, \lambda, \mu, \gamma \right)$ is additionally an M-stationary pair for CCOP. In order to show this, we consider the following cases:
\begin{itemize}
    \item [(i)] $\left\|\bar x\right\|_0 = s$. In this case, the assertion is trivially true, since the number of nonzero entries of $x$ must not exceed $s$. 
    \item[(ii)] $\left\|\bar x\right\|_0 < s$. We assume there exists an index $\hat i\in I_0\left(\bar x\right)\backslash I_0\left(x\right)$. Hence, $\bar x_{\hat i}=0$ and $x_{\hat i}\ne 0$. Hence, for the corresponding multiplier it holds $\bar \gamma_{\hat i}\ne 0$ since $\left(\bar x, \bar \lambda, \bar \mu, \bar \gamma \right)$ fulfils ND3. In a sufficiently small neighborhood, $\gamma_{\hat i}\ne 0$ holds. Since we assumed $\left( x, \lambda, \mu, \gamma \right)$ to be an M-stationary pair, it fulfils M3 from Definition \ref{def:M-stat}, i.\,e.
    $\gamma_i x_i=0$ for all $i=1, \ldots, n$. Consequently, $x_{\hat i}=0$, a contradiction. 
\end{itemize} 
Using the proven equality, $\left( x, \lambda, \mu, \gamma \right)$ is an M-stationary pair for CCOP if only if it holds, cf. Definition \ref{def:M-stat}:
\begin{itemize}
    \item[] M1: $D f(x) = \sum\limits_{p\in P}\lambda_p D h_p(x)+
    \sum\limits_{q \in Q} \mu_q D g_q( x)
    +\sum\limits_{i \in I_0\left(\bar x\right)} \gamma_i e_i$,
    \item[] M2: $\mu_q g_q(x)= 0$ and $\mu_q \geq 0$ for all $q \in Q$.
\end{itemize} 
Thus, $x$ is equivalently a Karusch-Kuhn-Tucker point for the following nonlinear program:
\begin{equation}
    \label{eq:locallysimp}
    \min_{x} \,\, f(x)\quad \mbox{s.\,t.} \quad x \in \bar M\end{equation}
with the feasible set given by equality and inequality constraints:
\[
    \bar M=\left\{x \in\R^{n} \left\vert\; h(x)=0, g(x)\ge 0, x_i=0, i \in I_0\left(\bar x\right) \right.\right\}.
\]
Hence, locally around $\left(\bar x, \bar \lambda, \bar \mu, \bar \gamma \right)$ CCOP becomes the easier optimization problem (\ref{eq:locallysimp}), which consists only of equality and inequality constraints.  Therefore, the standard result from Corollary 5.6 of \cite{kojima:1980} may be applied: $\bar x$ is strongly stable if and only if for all index subsets $Q_+(\bar x) \subset Q_* \subset Q_0(\bar x)$ the matrices $D^2 L(\bar x)\restriction_{\mathcal{T}_{\bar x} M_*(\bar x)}$ are nonsingular with the same determinant sign.

In order to prove the \emph{sufficiency part}, we assume that $\bar x$ is an M-stationary point and use the diffeomorphism $\Phi$ from (\ref{eq:stddiff}). We set $\bar f= f \circ \Phi^{-1}$ on the set $\left(\{0\}^{|P|} \times \H^{\left|Q_0(\bar x)\right|} \times \R^{n-k,s-k} \times \R^{k - m} \right)\cap V$ with $k=\|\bar x\|_0$. Recalling Remark \ref{rem:multipliers} we have at the origin:
\begin{itemize}
    \item[(i)] $\displaystyle \frac{\partial \bar f}{\partial y_{\left|P\right|+q}} \ge 0$ for $q \in Q_0\left(\bar x\right),$
    \item[(ii)]$\displaystyle \frac{\partial \bar f}
    {\partial y_{\ell+r}} = 0$ for $r=1,\ldots,n-\ell.$
\end{itemize}
Additionally we assume that $\bar x$ does not fulfil ND3, i.\,e. there exists $\bar i \in I_0(\bar x)$ such that $\gamma_{\bar i}=0$. Hence, we may assume without loss of generality that it holds:
\begin{itemize}
    \item [(iii)]$\displaystyle \frac{\partial \bar f}
    {\partial y_{\ell}} = 0. $
\end{itemize}
Moreover, we note that $s-k>0$, since ND3 is violated.
Next, we perturbate $\bar f$ by adding the term
\begin{equation}
    \label{eq:stabterm}
 \sum\limits_{q \in Q_0(\bar x)} \epsilon_{|P|+q} \cdot y_{|P|+q} + 
    \sum\limits_{i \in I_0(\bar x) \backslash \{k\}}  \epsilon_{m+i} \cdot y_{m+i} +\epsilon_{\ell} \cdot y_{\ell}^2+ 
    \sum\limits_{r \in R}  \epsilon_{\ell+r} \cdot y_{\ell+r}^2.
\end{equation}
We choose $\epsilon$'s in (\ref{eq:stabterm}) such that it holds for the perturbed function, which we again denote by $f$:
\begin{itemize}
    \item[(I)] $\displaystyle \frac{\partial f}{\partial y_{\left|P\right|+q}} > 0$ for $q\in Q_0\left(\bar x\right)$,
    \item[(II)]$\displaystyle \frac{\partial f}
    {\partial y_{\ell+r}} = 0$ for $r=0,\ldots,n-\ell$, and the matrix 
    $\displaystyle \left(\frac{\partial^2 f}{\partial y_{\ell+r_1} \partial y_{\ell+r_2}}\right)_{r_1,r_2 \in \left\{0,\ldots,n-\ell\right\}}$ is nonsingular,
    \item[(III)]$\displaystyle \frac{\partial f}
    {\partial y_{m+i}} \ne 0$ for $i \in I_0\left(\bar x\right)\backslash \{k\}.$
\end{itemize}
%Next, we approximate $f$ by a function $\bar f \in C^{\infty}$ which leaves its values and the derivatives of first and second order at the origin invariant. This can be done in some sufficiently small $C^2$-neighborhood of $f$ since $C^{\infty}$ lie $C^2$-dense in the $C^2$-functions. We denote the approximating function $\bar f$ by $f$ again.
Note, that this stabilizing step preserve the origin as an M-stationary point. 
We put 
\[
Y_{\ell-1}=\left(y_1,\ldots,y_{\ell-1}\right),Y_{n-\ell}=\left(y_{\ell+1},\ldots,y_n\right).\] 
Thus, $y=\left(Y_{\ell-1},y_{\ell},Y_{n-\ell}\right)$ for $y \in \left(\{0\}^{|P|} \times \H^{\left|Q_0(\bar x)\right|} \times \R^{n-k,s-k} \times \R^{k - m} \right)\cap V$.
It holds:
\[
\begin{array}{lcl}
f\left(Y_{\ell-1},y_{\ell},Y_{n-\ell}\right)&=&
\displaystyle \int_0^1 \frac{d}{dt}f\left(tY_{\ell-1},y_{\ell},Y_{n-\ell}\right)\diff t+f\left(0,y_{\ell},Y_{n-\ell}\right)\\ \\
&=&\displaystyle \sum\limits_{q \in Q_0(\bar x)} y_{|P|+q} \cdot d_{|P|+q}(y)+
    \sum\limits_{i \in I_0(\bar x)\backslash\{k\}}  y_{m+i} \cdot d_{m+i}(y)+ f\left(0,y_{\ell},Y_{n-\ell}\right),
\end{array}
\]
where
\[
\begin{array}{rcl}
   d_{|P|+q}(y)&=&\displaystyle \int_0^1 \frac{\partial f}{\partial y_{|P|+q}}\left(tY_{\ell-1},y_{\ell},Y_{n-\ell}\right)\diff t,\quad q \in Q_0(\bar x),
\\ \\
 d_{m+i}(y)&=&\displaystyle \int_0^1 \frac{\partial f}{\partial y_{m+i}}\left(tY_{\ell-1},y_{\ell},Y_{n-\ell}\right)\diff t,\quad i \in I_0(\bar x)\backslash\{k\}.   
\end{array}
\]
Note that 
 $d_{|P|+q} \in C^1$ for $q \in Q_0(\bar x)$,  and $d_{m+i} \in C^1$ for $i \in I_0(\bar x)\backslash\{k\}$.
Due to (II), we may apply the standard Morse Lemma on the $C^2$-function $f\left(0,y_{\ell},Y_{n-\ell}\right)$ without affecting the first $Y_{\ell-1}$ coordinates, see e.\,g. \cite{Jongen:2000}. The corresponding coordinate transformation is of class $C^1$. Denoting the transformed functions again by $f$ and $d_i$, we obtain
\[
f(y)= \sum\limits_{q \in Q_0(\bar x)} y_{|P|+q} \cdot d_{|P|+q}(y) +
    \sum\limits_{i \in I_0(\bar x)\backslash\{k\}}  y_{m+i} \cdot d_{m+i}(y) \pm y_{\ell}^2 + \sum\limits_{r \in R}\pm y_{\ell+r}^2.
\]
Conditions (I) and (III) provide that
\[
\begin{array}{rcl}
   d_{|P|+q}(0)&=&\displaystyle \frac{\partial f}{\partial y_{|P|+q}}\left(0\right) > 0,\quad q \in Q_0(\bar x),
\\ \\
 d_{m+i}(0)&=&\displaystyle \frac{\partial f}{\partial y_{m+i}}\left(0\right)\ne 0,\quad i \in I_0(\bar x)\backslash\{k\}.   
\end{array}
\]
Hence, we may take
\[
\begin{array}{l}
 y_{|P|+q} \cdot \left|d_{|P|+q}(y)\right|, \quad q \in Q_0(\bar x), \\
  y_{m+i} \cdot d_{m+i}(y), \quad  i \in I_0(\bar x)\backslash\{k\}, \\
  y_{\ell},\\
y_{\ell+r}, \quad r \in R
\end{array}
\]
as new local $C^1$-coordinates by a straightforward application of the inverse function theorem. We obtain in new coordinates locally around origin:
\[
  f(y)=\sum\limits_{q \in Q_0(\bar x)} y_{|P|+q} +
    \sum\limits_{i \in I_0(\bar x)\backslash\{k\}}  y_{m+i}  \pm y_{\ell}^2 + \sum\limits_{r \in R}\pm y_{\ell+r}^2.
\]
We perturb the resulting function for $\varepsilon >0$:
\[
f_{\varepsilon}(y)=\sum\limits_{q \in Q_0(\bar x)} y_{|P|+q} +
    \sum\limits_{i \in I_0(\bar x)\backslash\{k\}}  y_{m+i} \pm \left(y_{\ell} - \varepsilon\right)^2 + \sum\limits_{r \in R}\pm y_{\ell+r}^2.
\]
It is easy to see, that not only the origin, but also the point
\[
y_{\varepsilon}=(0,\ldots,\varepsilon,\ldots,0),
\]
where all coordinates but the $\ell$-th vanish, is M-stationary for $f_\varepsilon$. Especially, $y_\varepsilon$ is feasible since $s-k>0$ by assumption.
Overall, this shows that $f$ can be perturbed arbitrarily small to make a degenerate M-stationary point $\bar x$ bifurcate. Hence, $\bar x$ cannot be strongly stable. This concludes the proof.

\qed

Let us apply Theorem \ref{thm:char-ss} in order to characterize strongly stable local minimizers. Given an M-stationary point $\bar x \in M$  with mulitpliers $(\bar \lambda,\bar \mu, \bar \gamma)$ we set
\[
  M_+(\bar x) =\left\{
x \in \R^n\,\left\vert\,\begin{array}{l}
      h_p\left(x\right) =0, p \in P\\
g_q\left(x\right)=0,q \in Q_+(\bar x)\\
x_i=0, i \in I_0(\bar x)
\end{array} 
\right.\right\}.
\]
Obviously, $M_0(\bar x) \subset M_+(\bar x)$. In case that CC-LICQ holds at $\bar x$, the set $M_+(\bar x)$ is locally an  $C^2$-manifold of dimension $\left\|\bar x\right\|_0 - |P|-\left|Q_+(\bar x)\right|$. 
The tangent space of $M_+(\bar x)$ at $\bar x$ is thus given by
\[
\mathcal{T}_{\bar x} M_+(\bar x)=\left\{
\xi \in \R^n\,\left\vert\, \begin{array}{l}
   Dh_p(\bar x) \xi=0, p \in P \\
Dg_q(\bar x)\xi=0,q \in Q_+(\bar x)\\
\xi_i=0, i \in I_0(\bar x)
\end{array}
\right.\right\}.
\]

\begin{corollary}[Strongly stable minimizers]
\label{cor:ss-min}
Let $\bar x \in M$ be a strongly stable $M$-stationary point of CCOP satisfying CC-LICQ. Then, $\bar x$ is a local minimizer if and only if it fulfills ACC and the matrix $D^2 L(\bar x)\restriction_{\mathcal{T}_{\bar x} M_+(\bar x)}$ is positive definite.
\end{corollary}

\proof
First, we prove the \emph{if}-part. For that, let $\bar x$ fulfil ACC and let the matrix $D^2 L(\bar x)\restriction_{\mathcal{T}_{\bar x} M_+(\bar x)}$ be positive definite. Due to ACC, it holds $\left\|\bar x\right\|_0=s$.  
Moreover, for sufficiently close feasible points $x \in M$ it holds $I_0\left( x\right)\subset I_0\left(\bar x\right)$ due to continuity arguments and, thus, $I_0\left( x\right)= I_0\left(\bar x\right)$, since $\left\| x\right\|_0 \leq s$.
Hence, there exists a nonempty neighborhood $U_{\bar x}$ of $\bar x$ -- ensuring $x_i, i \in I_1(\bar x)$ do not vanish -- such that the feasible set locally becomes
 \[M\cap U_{\bar x} =\left\{x \in\R^n\cap U_{\bar x} \left\vert\; h(x)=0, g(x)\ge 0, x_i=0, i \in I_0\left(\bar x\right)\right.\right\}.
 \] Therefore, locally CCOP is a standard nonlinear program  which consists of equality and inequality constraints. We can then use the well-known sufficiency condition for its local minimizers, see e.\,g. Theorem 6 in \cite{mccormick:1967}. The latter states that the matrix $D^2 L(\bar x)\restriction_{\mathcal{T}_{\bar x} M_+(\bar x)}$ being positive definite is sufficient for $\bar x$ to be an (isolated) local minimizer.
 Next, we prove the \emph{only if}-part. For that, we consider two cases for the local minimizer $\bar x$ of CCOP:
 \begin{itemize}
     \item [(i)]ACC is fulfilled, i.\,e. $\left\|\bar x\right\|_0=s$. We follow the argumentation of the \emph{if}-part to conclude, that CCOP becomes locally a standard nonlinear program which consists of equality and inequality constraints. 
     %We can then use the known sufficiency condition for this, see e.\,g. Theorem 4 of \cite{mccormick:1967}, which leads to the matrix $D^2 L(\bar x)\restriction_{\mathcal{T}_{\bar x} M_0(\bar x)}$ being positive semi-definite as necessary condition. 
     Hence, $\bar x$ is an (isolated) local minimizer for the latter. Now, we apply Corollary 5.6 from \cite{kojima:1980} where strongly stable local minimizers were characterized in the context of nonlinear programming. In particular, it follows that $D^2 L(\bar x)\restriction_{\mathcal{T}_{\bar x} M_+(\bar x)}$ is positive definite.
     \item [(ii)] ACC is not fulfilled, i.\,e. $\left\|\bar x\right\|_0<s$.  We apply Proposition 2.1 from \cite{bucher:2018} to the local minimizer $\bar x$ of CCOP. The latter states that the corresponding Lagrange multipliers for the cardinality constraint vanish, i.\,e. $\gamma_i=0$ for all $i \in I_0(\bar x)$. This contradicts ND3 and, thus, due to Theorem \ref{thm:char-ss}, to the assumption that $\bar x$ is strongly stable. 
 \end{itemize}
\qed

Let us compare the cardinality-constrained second-order sufficient conditions (\ref{eq:cc-sosc32}) from \cite{bucher:2018} and (\ref{eq:cc-sosc321}) from \cite{pan:2017} with that characterizing strongly stable local minimizers $\bar x \in M$ from Corollary \ref{cor:ss-min}. If ACC holds at $\bar x$, then Lemma \ref{lem:crit-cone} is applicable and we obtain 
$\mathcal{C}_{\bar x} M \subset \mathcal{T}_{\bar x} M_+(\bar x)$.
Hence, the positive definiteness of the matrix $D^2 L(\bar x)\restriction_{\mathcal{T}_{\bar x} M_+(\bar x)}$ implies the cardinality-constrained second-order sufficient condition (\ref{eq:cc-sosc32}) from \cite{bucher:2018}. According to Proposition \ref{prop:opt-cond}, $\bar x$ is then a strict local minimizer. In view of (\ref{eq:under-acc}), we further recall that under ACC the cardinality-constrained sufficient optimality conditions (\ref{eq:cc-sosc32}) and (\ref{eq:cc-sosc321}) coincide. Hence, Proposition \ref{prop:opt-cond1} is also applicable and, thus, $\bar x$ is a strict local minimizer again.
However, the assumptions either in Proposition \ref{prop:opt-cond} or in Proposition \ref{prop:opt-cond1} (even if ACC is fulfilled) do not in general guarantee the strong stability of local minimizers. We illustrate this issue by means of the following Example \ref{ex:so-ss}.

\begin{example}[Sufficient optimality conditions and strong stability]
\label{ex:so-ss}
We consder the following CCOP with $P=\emptyset$, $Q=\{1,2\}$ and $n=3$, $s=2$:
\[
 \mbox{CCOP}: \quad \min_{x_1,x_2,x_3}\,\, f\left(x_1,x_2\right)=\left(x_1-1\right)^2 + 3\left(x_1-1\right)\left(x_2-1\right)+\left(x_2-1\right)^2 + x_3^2 \quad \mbox{s.\,t.}
\]
\[
 \quad 
   g_1\left(x_1,x_2,x_3\right)= x_1 -1 \geq 0, \quad
   g_2\left(x_1,x_2,x_3\right)= x_2 -1 \geq 0, \quad
   \left\|\left(x_1, x_2,x_3\right)\right\|_0 \leq 2.
\]
It is easy to see that the feasible point $\bar x=(1,1,0)$ solves $\mbox{CCOP}$. Since $\left\|\bar x\right\|_0=2$, ACC holds at $\bar x$. We have $Q_0(\bar x)=\{1,2\}$ and $I_0(\bar x)=\{3\}$, CC-LICQ is fulfilled at $\bar x$, and the corresponding multipliers vanish, i.\,e. $\bar \mu_1=\bar \mu_2=\bar \gamma_3=0$. We thus have $Q_+(\bar x)=\emptyset$ and $\mathcal{T}_{\bar x} M_+(\bar x)=\left\{\xi\in \R^3\,\left|\, \xi_3=0 \right.\right\}$. Let us compute the Hessian of the Lagrange function at $\bar x$:
\[
   D^2 L(\bar x) = \left(\begin{array}{ccc}
        2 &  3& 0\\
        3&  2& 0 \\
        0 & 0 & 2
   \end{array} \right).
\]
Note that its restriction on the tangent space of $M_+(\bar x)$ at $\bar x$ is indefinite:
\[
   D^2 L(\bar x) \restriction_{\mathcal{T}_{\bar x} M_+(\bar x)} = \left(\begin{array}{cc}
        2 &  3\\
        3&  2 \\
   \end{array} \right).
\]
Due to Corollary \ref{cor:ss-min}, $\bar x$ is not strongly stable. Now, we examine the cardinality-constrained second-order sufficient condition (\ref{eq:cc-sosc32}) from \cite{bucher:2018}. Lemma \ref{lem:crit-cone} provides the following representation of the critical cone of $M$ at $\bar x$:
\[
  \mathcal{C}_{\bar x} M=\left\{\xi\in \R^3\,\left|\, \xi_1 \geq 0, \xi_2 \geq 0, \xi_3=0 \right.\right\}.
\]
For all $\xi \in \mathcal{C}_{\bar x} M$ with $\xi\not=(0,0,0)$ it holds:
\[
    \xi^T \left(D^2f\left(\bar x\right) -
    \sum\limits_{q \in Q}\bar \mu_q D^2g_q\left(\bar x\right)\right) \xi = 2 \xi_1^2 + 6 \underbrace{\xi_1}_{\geq 0} \underbrace{\xi_2}_{\geq 0} + 2 \xi_2^2 > 0.
\]
We conclude that the cardinality-constrained second-order sufficient condition (\ref{eq:cc-sosc32}) from \cite{bucher:2018} is fulfilled at $\bar x$. Recall that under ACC the cardinality-constrained second-order sufficient condition (\ref{eq:cc-sosc321}) from \cite{pan:2017} is equivalent to the latter, thus, it also holds at $\bar x$. However, the minimizer $\bar x$ is not strongly stable for CCOP as we have seen before.  \qed
\end{example}

For the sake of completeness, we relate strong stability to nondegeneracy. The proof of the following Corollary \ref{cor:st-nd} is straightforward due to Theorem \ref{thm:char-ss}.  

\begin{corollary}[Stability and nondegeneracy]
\label{cor:st-nd}
Nondegenerate $M$-stationary points of CCOP are strongly stable. 
\end{corollary}

Since strongly stable M-stationary points are in particular locally unique, we immediately obtain the following result.  

\begin{corollary}[Local uniqueness and nondegeneracy]
\label{cor:lu-nd}
Nondegenerate $M$-stationary points of CCOP are locally unique. 
\end{corollary}

Note that Proposition \ref{prop:loc-uniq} and Corollary \ref{cor:lu-nd} are of independent interest and cannot be deduced one from each other. Next example shows that Corollary \ref{cor:lu-nd} may be well applied for M-stationary points, where the assumptions of Proposition \ref{prop:loc-uniq} fail to hold. This is in particular the case for M-stationary points of CCOP which are not local minimizers.

\begin{example}[Stability]
\label{ex:1}
We consider the following CCOP with $P=Q=\emptyset$ and $n=2$, $s=1$:
\[
 \mbox{CCOP}: \quad \min_{x_1,x_2}\,\, f\left(x_1,x_2\right)=(x_1-1)^2 + (x_2-1)^2 \quad \mbox{s.\,t.} \quad 
   \left\|\left(x_1, x_2\right)\right\|_0 \leq 1.
\]
It is easy to see that $\bar x=(0,0)$ is an M-stationary point, although not a local minimizer of $\mbox{CCOP}$. Nevertheless, it is nondegenerate. In fact, CC-LICQ is satisfied at $\bar x$, and ND1 holds. Since the cardinality constraint is not active: 
\[
   \left\|\bar x\right\|_0=0 < 1 =s,
\]
we have to check ND3:
\[
   \frac{\partial f}{\partial x_1} (\bar x) = -2, \quad \frac{\partial f}{\partial x_2} (\bar x) = -2.
\]
ND4 is trivially satisfied since $I_0(\bar x)=\{1,2\}$ and, hence, $T_{\bar x} M_0(\bar x)=\{(0,0)\}$. Overall, Corollary \ref{cor:lu-nd} applies for $\bar x$. Let us show that the 
cardinality-constrained second-order sufficient condition (\ref{eq:cc-sosc33}) from \cite{bucher:2018} is violated at $\bar x$. First, we have:
\[
\mathcal{C}_{\bar x} M=\left\{(\xi_1,\xi_2) \in \R^2 \,\left|\, \xi_1 + \xi_2 \geq 0, \left\|\left(\xi_1, \xi_2\right)\right\|_0 \leq 1\right. \right\}.
\]
But, for any $\xi \in \mathcal{C}_{\bar x} M$ with $\xi\not = (0,0)$ it holds:
\[
  \xi^T D^2f(\bar x) \xi = 0.
\]
We conclude that Proposition \ref{prop:loc-uniq} cannot be applied here.
\qed
\end{example}

\section{Global aspects}
\label{sec:glob}

We study the topological properties of lower level sets
\[
M_a=\left\{x\in M \left\vert f(x)\le a\right.\right\},
\]
where $a\in \R$ is varying. For that, we define intermediate sets for $a < b$:
\[
M^b_a=\left\{x\in M \left\vert a\le f(x) \le b\right. \right\}.
\]

\begin{assumption}
\label{as:proper}
The CCOP feasible set $M$ is compact and CC-LICQ is fulfilled at all points $x \in M$.
\end{assumption}

The deformation result for CCOP is based on the application of the standard Morse theory for nonlinear programming, see \cite{Jongen:2000}. For the sake of completeness we present the corresponding proof.

\begin{theorem}[Deformation for CCOP]
\label{thm:def}
Let Assumption \ref{as:proper} be fulfilled. If $M^b_a$ contains no M-stationary points for CCOP, then $M_a$ is homeomorphic to $M_b$.
\end{theorem}

\proof
For all $x\in M_a^b$ there exist due to CC-LICQ multipliers $\lambda_p(x),p\in P,\mu_q,q\in Q_0(x),\gamma_i(x),i\in I_0(x),\nu_r(x),r\in R(x)$ such that
\[
D f(x) = \sum\limits_{p\in P} \lambda_p D h_p( x)+
    \sum\limits_{q \in Q_0(x)} \mu_q D g_q(x)
    +\sum\limits_{i\in I_0(x)}^{n}  \gamma_i e_i + \sum\limits_{r\in R(x)}\nu_r\xi_r,
\]
where the vectors $\xi_r, r\in R(x)$ are chosen as in Lemma \ref{lem:local}.
Next, we set:
\[
\begin{array}{rcl}
A&=& \left\{x \in M_a^b \left\vert \mbox{there exists } r\in R(x) \mbox{ such that } \nu_r\ne 0\right. \right\},\\
B&=& \left\{x \in M_a^b \left\vert \mbox{there exists } q\in Q_0(x) \mbox{ such that } \mu_r< 0\right. \right\}.
\end{array}
\]
For $\bar x \in M_a^b$ we get $\bar x \in A \cup B$, since it is not M-stationary for CCOP.
The proof consists of a local argument an its globalization. First, we show the {\it local argument}, i.\,e. for each
$\bar x \in M_a^b$ there exist a neighborhood $U_{\bar x}$ of $\bar x, t_{\bar x}>0$, and a mapping
\[
\Psi_{\bar x}:\left\{
\begin{array}{rcl}
\left[0,t_{\bar x}\right) \times \left(M^b\cap U_{\bar x}\right)&\longrightarrow &M\\
\left(t,x\right)&\mapsto&\Psi_{\bar x}\left(t,x\right),
\end{array}\right.
\]
such that
\begin{itemize}
    \item [(i)] $\Psi_{\bar x}\left(t,M^b\cap U_{\bar x}\right)\subset M^{b-t}$ for all $t \in \left[0,t_{\bar x}\right)$,
    \item [(ii)] $\Psi_{\bar x}\left(t_1+t_2,\cdot\right)=\Psi_{\bar x}\left(t_1,\Psi_{\bar x}\left(t_2,\cdot\right)\right)$ for all $t_1,t_2 \in \left[0,t_{\bar x}\right)$ with $t_1+t_2 \in \left[0,t_{\bar x}\right)$,
    \item [(iii)] $\Psi_{\bar x}\left(\cdot,\cdot\right)$ is a $C^1$-flow corresponding to a $C^1$-vector field $F_{\bar x}$.
\end{itemize}
The level sets of $f$ are mapped locally onto the level sets of $f\circ \Phi^{-1}$, where $\Phi$ is the the diffeomorphism form Lemma \ref{lem:local}. We consider  $f\circ \Phi^{-1}$ and denote it by $f$ again.
Thus, we have $\bar x=0$ and $f$ is given on the feasible set
$\{0\}^{|P|} \times \H^{\left|Q_0(\bar x)\right|} \times \R^{n-k,s-k} \times \R^{k - m}$.\\
\begin{itemize}
    \item Case $\bar x \in A$. It follows from Remark \ref{rem:multipliers} that there exists $r\in R(\bar x)$ with $\displaystyle \frac{\partial f}{\partial x_r}\left(\bar x\right)\ne 0$. We define a local $C^1$-vector field $F_{\bar x}$ as
\[
F_{\bar x}(x_1,\ldots,x_r,\ldots,x_n)=\left(0,\ldots,\displaystyle-\frac{\partial f}{\partial x_r}\left(x\right)\cdot\left(\displaystyle\frac{\partial f}{\partial x_r}\left(x\right)\right)^{-2},\ldots,0\right)^T,
\]
which -- after respective inverse changes of local coordinates -- induces the flow $\Psi_{\bar x}$ fitting the local argument, see Theorem 2.7.6 from \cite{Jongen:2000}.
  \item Case $\bar x \in B$.
It follows from Remark \ref{rem:multipliers} that there exists $q\in Q_0(\bar x)$ with $\displaystyle \frac{\partial f}{\partial x_q}\left(\bar x\right)<0$. By means of a local $C^1$-coordinate transformation in the $q$-th coordinate on $\H$, leaving the other coordinates unchanged, we obtain locally for $f$
\[
f\left(x_1,\ldots,x_q,\ldots,x_n\right)=-x_q+f\left(x_1,\ldots,\bar x_q,\ldots,x_n\right).
\]
We define:
\[
F_{\bar x}(x_1,\ldots,x_q,\ldots,x_n)=\left(0,\ldots,1,\ldots,0\right)^T,
\]
which -- after respective inverse changes of local coordinates -- induces the flow $\Psi_{\bar x}$ fitting the local argument, see Theorem 3.3.25 from \cite{Jongen:2000}.
\end{itemize}

Next, we globalize the local argument.
For that, consider the open covering $\left\{U_{\bar x} \left\vert \bar x \in M_a^b \right.\right\}$. Since $M_a^b$ is compact, we get a finite subcovering $\left\{U_{\bar x_k} \left\vert \bar x_k \in M_a^b \right.\right\}$. Let $\left\{\phi_k\right\}$ be a $C^{\infty}$-partition of unity subordinate to this subcovering. We define the $C^1$-vector field
$F=\sum\limits_k \phi_kF_{x_k}$, which induces a flow on $\left\{U_{\bar x_k} \left\vert \bar x_k \in M_a^b \right.\right\}$.
%%%%%%%%%%%%
We define $r:[0,1]\times M^b\rightarrow M^b$ as
\[
r\left(\tau,x\right)=
\left\{
\begin{array}{ll}
x&\mbox{for }x \in M^a,\tau \in [0,1]\\
\Psi\left(\tau,x\right)&\mbox{for }x \in M_a^b,\tau \in [0,1].
\end{array}
\right.
\]
This mapping provides that $M_a$ is a strong deformation retract of $M^b$.
\qed

Let us now turn our attention to the topological changes of lower level sets when passing an M-stationary level. Traditionally, they are described by means of the so-called cell-attachment.
We first consider a special case of cell-attachment. For that, let $N^\epsilon$ denote the lower level set of a special linear function on $\H^w \times \R^{u,v}$, i.\,e.
\[
 N^{-\epsilon} = \left\{ \left(x,y\right) \in \H^w \times \R^{u,v} \,\left|\, \sum\limits_{i=1}^{w} x_i + \sum\limits_{j=1}^u y_j\leq -\epsilon \right. \right\},
\]
where $\epsilon \in \R$, and the integers $v < u$, and $w$ are nonnegative. 

\begin{lemma}[Normal Morse data]
\label{lem:cat}
For any $\epsilon > 0$ the set $N^\epsilon$ is homotopy-equivalent to $N^{-\epsilon}$ with $\binom{u-1}{v}$ cells of dimension $v$ attached. The latter cells are the $v$-dimensional simplices from the collection 
\[
 \left\{ \left. \mbox{conv} \left(e_j, j \in J\right)  \,\right|\, J \subset \{1, \ldots,p\}, 1 \in J, |J| = v+1 \right\}.
\]
\end{lemma}
\proof
 The lower level set $N^{-\epsilon}$ of a special linear function on $\H^w \times \R^{u,v}$ is given by
\[
 N^{-\epsilon} = \left\{ \left(x,y\right) \in \H^w \times \R^{u,v} \,\left|\, \sum\limits_{i=1}^{w} x_i + \sum\limits_{j=1}^u y_j\leq -\epsilon \right. \right\}
\]
with $\epsilon >0$. This is homotopy-equivalent to the set
\[
 \bar N^{-\epsilon} = \left\{ y \in \R^{u,v} \,\left|\, \sum\limits_{i=1}^{u} y_i \leq -\epsilon \right. \right\}, 
\]
using the homotopy
\[
\left((x,y),t\right)) \mapsto \left((1-t) \cdot x,y\right), \quad t \in [0,1].
\]
We note that $\bar N^{-\epsilon}$ is a lower level set of the special problem of sparsity constrained nonlinear optimization
\begin{equation}
\label{eq:normalold}
\min_{y \in \R^u} \,\sum\limits_{i=1}^{u} y_i \, \quad \mbox{s.\,t.} \quad \NNorm{y}{0} \le v.
\end{equation}
The Cell-Attachmemt for (\ref{eq:normalold}) was examined in \cite{laemmel:2019}. According to the latter,
$ \bar N^{-\epsilon}$ is homotopy-equivalent to $ \bar N^{\epsilon}$ with $\binom{u-1}{v}$ cells of dimension $v$ attached. The set $ \bar N^{\epsilon}$ is homotopy-equivalent to the set $N^{\epsilon}$ by using the same homotopy as above. Thus, the assertion follows immediately.

 \qed

\begin{theorem}[Cell-Attachment for CCOP]
\label{thm:cell-a}
Let Assumption \ref{as:proper} be fulfilled. Suppose that $M^b_a$ contains exactly one nondegenerate M-stationary point $\bar x$ with $\left\|\bar x\right\|_0=k$ and the M-index equal to $s-k+QI$. 
If $a<f\left(\bar x \right) <b$,
then $M^b$ is homotopy-equivalent to $M^a$ with $\binom{n-k-1}{s-k}$ cells of dimension $s-k+QI$ attached, namely:
\[
 \bigcup_{\scriptsize \begin{array}{c}
        J \subset \left\{1, \ldots,n-k\right\} \\
        1 \in J, |J| = s-k +1
   \end{array}} \mbox{conv} \left(e_j, j \in J\right)\times [0,1]^{QI}.
\]
\end{theorem}

\proof Theorem \ref{thm:def} allows deformations up to an arbitrarily small neighborhood of the M-stationary point $\bar x$. In such a neighborhood, we may assume without loss of generality that $\bar x=0$ and $f$ has the following form as from Theorem \ref{thm:morse}:
\begin{equation}
\label{eq:n1}
    f(x)= f(\bar x) +
    \sum\limits_{q \in Q_0(\bar x)} x_{|P|+q} +
    \sum\limits_{i \in I_0(\bar x)}  x_{m+i} + \sum\limits_{r \in R}\pm x_{\ell+r}^2,  
\end{equation}
where $x \in \{0\}^{|P|} \times \H^{\left|Q_0(\bar x)\right|} \times \R^{n-k,s-k} \times \R^{k - m}$. Moreover, there are exactly $QI$ negative squares in (\ref{eq:n1}).
In terms of \cite{goresky:1988} the set 
$\{0\}^{|P|} \times \H^{\left|Q_0(\bar x)\right|} \times \R^{n-k,s-k} \times \R^{k - m}$ can be interpreted as the product of the tangential part $\{0\}^{|P|}  \times \R^{k - m}$ and the normal part $\H^{\left|Q_0(\bar x)\right|} \times  \R^{n-k,s-k}$. The cell-attachment along the tangential part is standard.  Analogously to the case of nonlinear programming, one $QI$-dimensional cell has to be attached on $\{0\}^{|P|} \times \R^{k - m}$.  The cell-attachment along the normal part is more involved. Due to Lemma \ref{lem:cat}, we need to attach $\binom{n-k-1}{s-k}$ cells on $\H^{\left|Q_0(\bar x)\right|} \times\R^{n-k,s-k}$, each of dimension $s-k$. Finally, we apply Theorem 3.7 from Part I in \cite{goresky:1988}, which says that the local Morse data is the product of tangential and normal Morse data. Hence, the dimensions of the attached cells add together. Here, we have then to attach $\binom{n-k-1}{s-k}$ cells on $\{0\}^{|P|} \times \H^{\left|Q_0(\bar x)\right|} \times \R^{n-k,s-k} \times \R^{k - m}$, each of dimension $s-k+QI$. \qed

Let us present a global interpretation of Theorems \ref{thm:def} and \ref{thm:cell-a} by deducing a mountain pass result for CCOP.

\begin{remark}[Mountain Pass]
\label{rem:m-pass}
Let Assumption \ref{as:proper} be fulfilled and additionally the CCOP feasible set $M$ be connected. Then, it holds:
%the number of connected components of $M^a$ is equal $0$ for $a=a_1$ small enough and equal to $1$ for $a=a_2$ big enough. The change in the number of connected component of the lower level sets in between, i.\,e. of $M^b, a_1\le b\le a_2$, can be estimated in terms of 
\[
r_1+(n-s)r_2\ge r-1,\]
where $r$ is the number of local minimizers for CCOP, $r_1$ is the number of its saddle points $\bar x$ with M-index equal to 1 and active sparsity constraint, i.\,e. $\left\|\bar x\right\|_0=s$, and $r_2$ is the number of saddle points $\bar x$ with M-index equal to 1 and inactive sparsity constraint, i.\,e. $\left\|\bar x\right\|_0=s-1$. This is due to the following implications of Theorem \ref{thm:cell-a}:
\begin{itemize}
    \item [(1)] at most one component disappears when passing a saddle point with M-index equal to 1 and active sparsity constraint,
    \item[(2)] at most $n-s$ components disappear when passing a saddle points with M-index equal to 1 and inactive sparsity constraint,
    \item[(3)] a new component is created when passing a local minimizer, and
    \item[(4)] no change in the number of components happens when passing any other point.
\end{itemize}
Moreover, for a sufficiently small level the corresponding lower level set is empty. For a sufficiently large level the corresponding lower level set is connected. \qed
\end{remark}

\section*{Acknowledgment}
The authors would like to thank Hubertus Th. Jongen for fruitful discussions.

\bibliographystyle{apalike}
\bibliography{lit.bib}

\end{document}